\documentclass[preprint]{elsarticle}

\usepackage[pdftex,dvipsnames]{xcolor}  
\usepackage[prependcaption,textsize=tiny]{todonotes}  
\usepackage{url}
\usepackage{graphicx}
\usepackage{subfigure}
\usepackage{amsmath, amssymb, amsthm,scalerel}
\usepackage{cases}
\usepackage{epstopdf}

\usepackage{colortbl}
\usepackage{tikz}
\usetikzlibrary{matrix, calc, patterns}
\usepackage{listings}

\newcommand{\acrfont}{\fontfamily{qcr}\selectfont }
\newtheoremstyle{example_style}
{4pt}
{4pt}
{}
{}
{\bfseries}
{:}
{.5em}
{}
\theoremstyle{example_style}
\newtheorem{example}{Example}
\newcommand{\bu}{{\boldsymbol u}}
\newcommand{\bx}{{\boldsymbol x}}
\DeclareMathOperator*{\Bigcdot}{\scalerel*{\cdot}{\bigodot}}
\usepackage{hyperref}
\usepackage[all]{hypcap}

\begin{document}

%
\begin{frontmatter}
\title{ParMooN -- a modernized program package based on mapped finite elements}
\author[Volker1]{Ulrich Wilbrandt}
\ead{ulrich.wilbrandt@wias-berlin.de}
\author[Volker1]{Clemens Bartsch}
\ead{clemens.bartsch@wias-berlin.de}
\author[Volker1]{Naveed Ahmed}
\ead{naveed.ahmed@wias-berlin.de}
\author[Volker1,Najib]{Najib Alia}
\ead{najib.alia@wias-berlin.de}
\author[Volker1,Felix]{Felix Anker}
\ead{felix.anker@wias-berlin.de}
\author[Volker1]{Laura Blank}
\ead{laura.blank@wias-berlin.de}
\author[Volker1]{Alfonso Caiazzo}
\ead{alfonso.caiazzo@wias-berlin.de}
\author[Sashi,sashi]{Sashikumaar Ganesan}
\ead{sashi@cds.iisc.ac.in}
\author[Volker1]{Swetlana Giere}
\ead{swetlana.giere@wias-berlin.de}
\author[Gunar]{Gunar Matthies}
\ead{gunar.matthies@tu-dresden.de}
\author[Sashi]{Raviteja Meesala}
\ead{raviteja@cmg.cds.iisc.ac.in}
\author[Sashi]{Abdus Shamim}
\ead{shamim@cmg.cds.iisc.ac.in} 
\author[Sashi]{Jagannath Venkatesan}
\ead{jagan@cmg.cds.iisc.ac.in} 
\author[Volker1,Volker2,Volker3]{Volker John\corref{cor}}
\ead{volker.john@wias-berlin.de}

\cortext[cor]{Corresponding author.}
\address[Volker1]{Weierstrass Institute for Applied Analysis and Stochastics,
Leibniz Institute in Forschungsverbund Berlin e. V. (WIAS), Mohrenstr. 39, 10117 Berlin, Germany}
\fntext[Najib]{The work of Najib Alia has been supported by a funding from the 
European Union's Horizon 2020 research and innovation programme under the Marie Sk{\l}odowska-Curie grant agreement No. 675715 (MIMESIS).}

\fntext[Felix]{The work of Felix Anker has been 
supported by grant Jo329/10-2
  within the DFG priority programme 1679: Dynamic simulation of interconnected 
  solids processes.}
\address[Sashi]{Department of Computational and Data Sciences,  
Indian Institute of Science,\\ Bangalore - 560012, India}
\fntext[sashi]{The work of Sashikumaar Ganesan has  partially been supported by the Naval Research Board, DRDO, India through the grant NRB/4003/PG/368.}
\address[Gunar]{
Department of Mathematics,
Institute of Numerical Mathematics, TU Dresden,\\
01062 Dresden
Germany}
\address[Volker2]{Free University of Berlin,
Department of Mathematics and Computer Science,
Arnimallee 6, 14195 Berlin, Germany}
\fntext[Volker3]{The work of Volker John has  partially been
supported by grant Jo329/10-2
  within the DFG priority programme 1679: Dynamic simulation of interconnected 
  solids processes.\\[1ex]
  \textcopyright~2017. This manuscript version is made available under the CC-BY-NC-ND 4.0 license http://creativecommons.org/licenses/by-nc-nd/4.0/}


\begin{abstract}
{\sc ParMooN} is a program package for the numerical solution of elliptic and 
parabolic partial differential equations. It inherits the distinct features 
of its predecessor {\sc MooNMD} \cite{JM04}: strict decoupling of geometry and finite element 
spaces, implementation of mapped finite elements as their definition can be 
found in textbooks, and a geometric multigrid preconditioner with the option 
to use different finite element spaces on different levels of the multigrid 
hierarchy. After having presented some thoughts about in-house research codes, 
this paper focuses on aspects of the parallelization for a distributed memory environment, which is the main novelty 
of {\sc ParMooN}. 
Numerical studies, performed on compute servers, assess the efficiency of the 
parallelized geometric multigrid preconditioner in comparison with some parallel
solvers that are available in the library {\sc PETSc}. The results of these studies give a first
indication whether the cumbersome implementation of the parallelized geometric
multigrid method was worthwhile or not. 
\end{abstract}


\end{frontmatter}

\section{Introduction}

{\sc MooNMD}, a  C++ program package for the numerical solution of 
elliptic and parabolic partial differential equations based on mapped finite elements, is described in 
\cite{JM04}.
A modernized version of this package, called  {\sc ParMooN}, has  recently been developed 
to be used as a research code in the future. 

The core of  {\sc MooNMD} was designed more than 15 years ago and this code has been successfully used 
in many scientific studies. There are almost 90 research articles citing {\sc MooNMD} via \cite{JM04}, see \cite{scopus}. 
Recent advances in computing hardware and language standards 
necessitate a re-design and re-implementation of some of the core routines. 
With the new core and the new features, the code 
was renamed to {\sc ParMooN} (Parallel Mathematics and object-oriented Numerics).

The general aims of this paper are 
to report on the development of the existing research code
towards the new package {\sc ParMooN} in order to accomplish
the desired features of an in-house research code  that will be formulated in Section~\ref{sec:research_codes} and to assess the parallelized 
geometric multigrid method by comparing it with solvers that are available in an external library. 
The original code possesses some distinct features that should be transferred to {\sc ParMooN},
like the strict decoupling of 
geometry and finite element spaces, the implementation of mapped finite elements as their definition
can be found in textbooks, and a multiple discretization multilevel (MDML) preconditioner.
This paper focuses on the most relevant aspect concerning the development of  {\sc ParMooN},
namely the distributed memory parallelization. In particular, 
the technically most cumbersome part, the parallelization of the geometric multigrid, 
is discussed.

A main contribution of this paper is a first assessment of the resulting parallel geometric multigrid method 
in comparison with parallel solvers for linear systems of equations that can be called from 
the library {\sc PETSc}.
The numerical studies were performed on compute servers as the available in-house hardware.
We think that this assessment is also of interest for other 
groups who develop their own codes in order to get an impression whether it is 
worthwhile to implement a parallelized geometric multigrid method 
or not. Two main problem classes supported in {\sc ParMooN} are considered in the 
numerical studies: scalar convection-diffusion-reaction equations and 
the incompressible Navier--Stokes equations.

The paper is organized as follows. Section~\ref{sec:research_codes} contains an exposition of our 
thoughts about in-house research codes, in particular about their advantages and their goals. Mapped finite 
elements, as they are used in {\sc MooNMD}/{\sc ParMooN}, are described in Section~\ref{sec:mapped_fe}.
Section~\ref{sec:par_data} presents main aspects of the parallelization and the parallelization of the 
geometric multigrid method is briefly discussed in Section~\ref{sec:par_mg}. Numerical studies that 
compare this method with solvers available in  {\sc PETSc} are presented in Section~\ref{sec:numres}. The paper concludes with a summary.

\section{Some Considerations about In-House Research Codes}\label{sec:research_codes}

Nowadays, several academic software packages for 
solving partial differential equations exist in the research community. They are 
usually developed and supported by research groups for whom
software development is one of the main scientific tasks. 
Such software packages include, among others, 
{\sc deal.II} \cite{dealII16}, {\sc FEniCS} \cite{fenics15}, {\sc DUNE} \cite{DKNO10,BBD+16}, 
{\sc OpenFOAM} \cite{openfoam}, or {\sc FreeFem++} \cite{FreeFem}. 
These packages have advanced functionality and support 
features 
like adaptive mesh refinement, parallelism, etc.

Naturally, a research code developed in-house possesses less functionality than these
large packages.
In view of their availability, the following questions arise:
Why is it worth to develop an own research code? In particular, 
is it worth to develop a code within a research group that focuses primarily on
numerical analysis?
 In the following, some arguments, mainly based on our own experience, are presented. 

{\em In-depth knowledge of details of the software.}
The first key aspect of working with a code developed and maintained within the research group
is the detailed knowledge of the software structure.
In fact, applied mathematicians often work at the development of numerical methods. 
These methods have to 
be implemented, assessed, and compared with popular state-of-the-art methods for the same problem. 
A meaningful assessment requires the usage of the methods in the same code. In this respect, 
it is important
to have access to a code where one knows and can control every detail.

For brevity, just one example will be mentioned to show the importance of knowing the details of a software package.
This example concerns the clarification of 
appropriate interface conditions in subdomain iterations for the Stokes--Darcy problem, see \cite{CJW14}.
Standard Neumann interface conditions can be used only for viscosity and permeability coefficients 
that are unrealistically large. For realistic coefficients, appropriate Robin boundary conditions
have to be used. The implementation of the Robin interface conditions was performed in a straightforward way in {\sc MooNMD}/{\sc ParMooN}. 

{\em Flexibility.}
Further advantages of an own research code are the possibility of controlling its core parts and flexibility. In particular, for our research it is 
very important that the code supports the use of different discretization strategies.
As an example,  {\sc MooNMD} 
was designed for finite element methods. But for the investigation of discretizations of time-dependent
convection-diffusion equations in \cite{JN12}, finite difference methods were implemented
as well. Because these methods performed very well, they were later used 
in the context of simulating 
population balance systems defined in tensor-product domains, e.g., see \cite{SBTJ14}.

{\em Testing of numerical methods.}
{\sc MooNMD} was used in the definition of benchmark problems in \cite{HTK09}.
The list of examples could be extended. 
In addition, a number of numerical methods have been developed 
and implemented in our research code which turned out to be not (yet) competitive, like the optimization of stabilization
parameters in SUPG methods in \cite{JKS11}. 
Having  a known and flexible research code 
at disposal allows to test and support 
methods that, at the time of the implementation, have not been benchmarked in detail.

Certainly, also the large packages mentioned above allow the implementation of
different methods and discretization strategies. However, we think that a successful implementation often
requires a very close interaction with core developers of the packages. 
This effort might not be feasible for both the user and the developer. 
Therefore, an own code might reduce
the time from the development to the assessment of numerical methods. 

{\em Benefits for students.}
A further aspect, related to the interaction of the core development team, 
concerns the students and the PhD students who are involved in the 
development and in 
the usage of the code.
Since the core developers of an own research code are readily available and they are experts in the focused research topics of 
the group, these students can be 
supported efficiently. In addition, students working at the code stated 
several positive effects: the work at details of the implementation 
facilitates the insight into the methods and algorithms, 
which is important for analyzing their properties, and it enhances the 
skills in software design and management. 

Of course, incorporating students into 
code development requires that there is an easy use of the code and an easy access to basic routines, 
such as, in the case of a finite element solver, assembling of matrices and solving linear systems of equations. 
This issue touches already the next question: What should be expected from an own research code? 

{\em Easy usage.}
In order to support students starting to work with the code, an easy installation and 
basic testing setup are essential.
There are even successful attempts for designing complex codes that can 
be used for teaching students in basic courses on numerical methods for partial 
differential equations, like the so-called computational laboratory 
for Investigating Incompressible Flow Problems ({\sc IFISS}), see \cite{ERS14}, which uses {\sc Matlab}, and the open source software {\sc FreeFem++} \cite{FreeFem}, which is based on an own language.

{\em Modularity.}
The code should be modular. In particular, there should be a general core and 
individual projects are attached to this core. Of course, the projects use routines from the core. But 
using
an own code, it seems to be easier than with a large package, which  is developed somewhere else, to
incorporate contributions from the projects into the core. 

{\em Stability.}
With respect to the required stability of the code,  there are, in our opinion, no fundamental differences
between own research codes and large packages. 

{\em Efficiency.}
However, there are different expectations with respect to the efficiency. A research code should be 
flexible in many respects, since its main tasks include supporting the development of numerical methods 
and results from numerical analysis. For instance, in the code {\sc MooNMD}/{\sc ParMooN}, the concept 
of mapped finite elements is implemented, see Section~\ref{sec:mapped_fe}. In this way, the code supports
currently around 170 finite elements in two dimensions and 75 finite elements in three dimensions.
Consequently, all routines are implemented for the general situation. For certain finite elements, this 
might be less efficient than 
using tailored routines. 
However,
also for a research code, efficiency is a key property that should not be neglected.
For instance, 
the simulation of standard academic benchmark problems for turbulent incompressible flows requires 
the computation of large time intervals to collect temporal averages of statistics of interest. 
In our opinion, an own research code should be reasonably efficient on the available in-house hardware, 
which, in our case, are usually laptops, compute servers, or small clusters.

\section{Mapped Finite Element Spaces}\label{sec:mapped_fe}

The implementation of finite element methods
in {\sc MooNMD}/{\sc ParMooN} is based on a rather 
abstract definition of a finite element space and on the mapping of each mesh cell to a 
reference cell. 

Let $\Omega \subset \mathbb R^d$, $d\in \{2,3\}$, be a bounded domain and let $\mathcal T^h$ be an admissible
triangulation of $\Omega$ consisting of compact, simply connected mesh cells. For each mesh cell $K$, 
a local finite element space $P(K)\subset C^s(K)$, $s\ge  0$, 
is given by some finite-dimensional space of functions
spanned by a basis 
$\{\phi_{K,i}\}_{i=1}^{N_K}$. Furthermore, a set of local linear functionals  
$\{\Phi_{K,i}\}_{i=1}^{N_K}$ is given. The space $P(K)$ is unisolvent with respect to the 
functionals. Often, a so-called local basis is chosen, i.e., a basis that satisfies 
$\Phi_{K,i}(\phi_{K,j}) = \delta_{ij}$ for  $i,j=1,\ldots,N_K$. The local linear functionals might 
be values of the functions or their derivatives in certain points, integrals on $K$ or on faces of
$K$.

Let $\Phi_1,\ldots,\Phi_N\,
: \, \left\{v\in L^\infty(\Omega)\ : \ v|_K \in P(K)\right\} \to \mathbb R$ be given continuous linear functionals, where $v|_K \in P(K)$ has to be understood in the way that $v$ is a polynomial in the 
interior of $K$ and it is extended continuously to the boundary of $K$. 
The restriction of each functional to $C^s(K)$ defines the set of local functionals.
The union of all mesh cells $K_j$,
for which there is a $p\in P(K_j)$ with $\Phi_i(p) \neq 0$, will be denoted by
$\omega_i$. Now, 
the global finite element space is defined as follows. 
A function $v(\bx)$ defined on $\Omega$ with 
  $v|_K \in P(K)$ for all $K \in \mathcal T^h$ is called continuous with respect
to the functional $\Phi_i$ if 
$\Phi_i(v|_{K_1}) = \Phi_i(v|_{K_2})$ for all $K_1, K_2 \in \omega_i$.
The space
\begin{eqnarray*}
S &=& \Big\{ v\in L^\infty(\Omega)\ : \ v|_K \in P(K) \text{ and }
v \text{ is continuous with respect to }\\
&& \hspace*{0.7em} \Phi_i, i=1,\ldots,N\Big\}
\end{eqnarray*}
is called finite element space.
The global basis 
$\{\phi_j\}_{j=1}^N$ of $S$ is defined by the 
conditions
$\phi_j\in S$ with $\Phi_i(\phi_j) = \delta_{ij}$ for  $i,j = 1,\ldots,N$.

Using this definition for the implementation of a finite element space requires
\begin{list}{}{\itemsep0.0ex\parsep0.1ex\topsep0.0ex}
\item[1.)] the definition of 
the local basis and linear functionals for each $K$,
\item[2.)] the implementation of a method that assures continuity  with respect to the 
functionals stated in the definition of~$S$.
\label{page:point2}
\end{list}

The first requirement can be achieved in two different ways, via a mapped or an unmapped
implementation.
In the unmapped approach, the local basis and linear functionals are defined directly on $K$. 
In contrast, mapped finite elements are closely connected to a standard way of
analyzing finite element discretizations. This analysis consists of
three steps:
\begin{list}{}{\itemsep0.0ex\parsep0.1ex\topsep0.0ex}
\item[$\bullet$] Map an arbitrary mesh cell $K$ to a compact reference mesh cell
  $\hat K$.
\item[$\bullet$] Prove the desired properties on $\hat K$, which is the core of the analysis.
\item[$\bullet$] Map the reference mesh cell $\hat K$ back to $K$ to get the
  final result.
\end{list}
Hence, this approach has two main features:
\begin{list}{}{\itemsep0.0ex\parsep0.1ex\topsep0.0ex}
\item[$\bullet$] All considerations have to be done on $\hat K$ only.
\item[$\bullet$] Information about neighbor mesh cells of $K$ is neither available nor needed.
\end{list}
Mapped and unmapped finite element methods possess the same analytical
properties if the reference map $F_K \, :\, \hat K \to
K$ is affine for every mesh cell $K$ of the given
triangulation, e.g., compare \cite[Chap.~2.3]{Cia78}. In the case of non-affine maps, occurring, e.g., 
for a triangulation consisting of arbitrary quadrilateral or hexahedral
mesh cells, mapped and unmapped finite element spaces might be
different. 
In {\sc MooNMD}/ {\sc ParMooN}, the concept of mapped finite elements is implemented in the 
following way. Reference mesh
cells are the unit simplices, e.g., in two dimensions with the vertices $(0,0)$, $(1,0)$, $(0,1)$, 
and the unit cubes $\hat K = [-1,1]^d$. Affine maps are available for all reference mesh cells. 
To account for arbitrary quadrilaterals and hexahedra, $d$-linear maps are also implemented for 
the unit cubes. Based on the different reference cells, local spaces on $\hat K$, linear functionals, and reference maps, 
{\sc MooNMD}/{\sc ParMooN} currently supports about 170 finite elements in two dimensions and 75 finite elements in 
three dimensions.

The use of mapped finite element spaces  essentially requires  the implementation of 
finite elements on the reference cells. The quadrature rules for numerical integration
have to be implemented only on these cells, since the integrals on physical cells are
transformed to integrals on the reference cells. Note that the same strategy works also for
the handling of cell faces, which are mapped onto lower-dimensional reference cells by corresponding reference maps.

Concerning requirement 2.) above, a finite element space $S$ is represented by 
a map $\mathcal F$ called d.o.f.-manager which maps local, i.\,e., within a 
cell, indices  of degrees of freedom (d.o.f.) to global ones. For this purpose define 
$M(K)$ to be the set of local d.o.f.s denoted by $(K, i)$ on the cell $K$. Then 
define the set of all local degrees of freedom
\begin{equation*}
 M := \bigcup_{K\in \mathcal T^h} M(K).
\end{equation*}
The local-to-global map $\mathcal F$ now surjectively maps $M$ to 
$\{1,\dots,N\}$ such that $\mathcal F((K,i)) = \mathcal F((K',j))$ whenever
the local degrees of freedom $(K,i)$ and $(K',j)$ belong to the same 
global degree of freedom. The number 
$N \le |M|$ is then the number of global degrees of freedom. In other words, 
$\mathcal F$ describes a partition $P$ of the set $M$ (i.\,e., an equivalence
relation) together with a global numbering. Computing such a map $\mathcal F$ 
is done via Algorithm~\ref{alg:local_to_global_dof_map} in {\sc ParMooN}. 
{
\renewcommand\lstlistingname{Algorithm} 
\begin{lstlisting}[
float,
caption={Computation of $\mathcal F$ that maps local degrees of freedom 
         to global ones. The input is a mesh $\mathcal T^h$ whose cells $K$ are 
         ordered by increasing integers $id(K)$.},
label=alg:local_to_global_dof_map,
xleftmargin=4em,
language=Python,
numbers=left,
mathescape=true,
numberstyle=\tiny,
basicstyle=\small,          % print whole listing small
keywordstyle=\color{black}\bfseries, % bold black keywords
identifierstyle=,           % nothing happens
commentstyle=\color{white}, % white comments
stringstyle=\ttfamily,      % typewriter type for strings
showstringspaces=false,     % no special string spaces
escapeinside={(*@}{@*)}
]
$P \leftarrow \text{ finest partition of }M$
for $\text{all mesh cells } K\in \mathcal T^h$
   $\text{determine } M(K)$
   for $\text{all neighbors } K' \text{ of } K \text{ with } id(K) < id(K')$
      $\text{determine } M(K')$
      $\text{find local partition } P_{\mathrm{loc}} \text{ of } M(K)\cup M(K')$  (*@\label{alg_line:find_partitioning}@*)
      $P \leftarrow \text{finest partition coarser than } P \text{ and } P_{\mathrm{loc}}$(*@\label{alg_line:update_partitioning}@*)
$\text{assign increasing integers to each subset in the partition }P$
\end{lstlisting}
}

Initially the partition is set to be the finest partition of the set $M$, 
i.\,e., it consists of disjoint single-element subsets, so that no two degrees 
of freedom in $M$ are identified yet. Step~\ref{alg_line:find_partitioning} in 
Algorithm~\ref{alg:local_to_global_dof_map} consists of finding identical 
degrees of freedom in two neighboring cells $K$ and $K'$, i\,e., a local 
partition $P_{\mathrm{loc}}$ of $M(K)\cup M(K')$. This  step is done using
the information of positions of degrees of freedom in the reference cell and 
tailored mapper classes. Step~\ref{alg_line:update_partitioning} updates the 
partition $P$ such that all sets in the previous set $P$ and in $P_{\mathrm{loc}}$
are contained in one set in the updated partition $P$. Finally, to each set of 
the resulting partition of $M$, a unique integer is assigned that will serve as 
(the index of) the global degree of freedom. 
Steps~\ref{alg_line:find_partitioning} and \ref{alg_line:update_partitioning} 
are explained in detail in the following example.

\begin{example} \label{ex:local_to_global_dof_map}
 Consider a $2\times2$ mesh $\mathcal T^h$ consisting of four cells $A$, $B$, 
 $C$, and $D$, each associated with a $Q_1$-finite element, see 
 Figure~\ref{fig:local_to_global_dof_mesh}. 
 \begin{figure}
 \begin{center}
\newcommand{\twoByTwoGridWithCellNames}{
 \coordinate (bottom_left) at (0,0);
 \coordinate (bottom_right) at (2,0);
 \coordinate (top_left) at (0,2);
 \coordinate (top_right) at (2,2);
 \coordinate (bottom) at ($(bottom_left)!0.5!(bottom_right)$);
 \coordinate (top) at ($(top_left)!0.5!(top_right)$);
 \coordinate (left) at ($(bottom_left)!0.5!(top_left)$);
 \coordinate (right) at ($(top_right)!0.5!(bottom_right)$);
 \coordinate (center) at ($(left)!0.5!(right)$);
 %
 \draw[-] (top_left) -- (top_right);+
 \draw[-] (left) -- (right);
 \draw[-] (bottom_left) -- (bottom_right);
 \draw[-] (top_left) -- (bottom_left);
 \draw[-] (top) -- (bottom);
 \draw[-] (top_right) -- (bottom_right);
 \draw ($(left)!0.5!(bottom)$) node {A};
 \draw ($(right)!0.5!(bottom)$) node {B};
 \draw ($(left)!0.5!(top)$) node {C};
 \draw ($(right)!0.5!(top)$) node {D};
}
\newcommand{\markCellA}[4]{
 \draw (bottom_left) node[above right] {#1};
 \draw (bottom) node[above left] {#2};
 \draw (left) node[below right] {#3};
 \draw (center) node[below left] {#4};
}
\newcommand{\markCellB}[4]{
 \draw (bottom) node[above right] {#1};
 \draw (bottom_right) node[above left] {#2};
 \draw (center) node[below right] {#3};
 \draw (right) node[below left] {#4};
}
\newcommand{\markCellC}[4]{
 \draw (left) node[above right] {#1};
 \draw (center) node[above left] {#2};
 \draw (top_left) node[below right] {#3};
 \draw (top) node[below left] {#4};
}
\newcommand{\markCellD}[4]{
 \draw (center) node[above right] {#1};
 \draw (right) node[above left] {#2};
 \draw (top) node[below right] {#3};
 \draw (top_right) node[below left] {#4};
}
\begin{tikzpicture}[ scale=2 ]
\begin{scope}
 \twoByTwoGridWithCellNames
 \markCellA{1}{2}{3}{4}
 \markCellB{1}{2}{3}{4}
 \markCellC{1}{2}{3}{4}
 \markCellD{1}{2}{3}{4}
\end{scope}

\draw(2.2,1) edge[->,bend right=-45] node[auto]{$\mathcal F$}  (2.8,1);

\begin{scope}[shift={(3,0)}]
 \twoByTwoGridWithCellNames
 \markCellA{1}{2}{3}{4}
 \markCellB{2}{5}{4}{6}
 \markCellC{3}{4}{7}{8}
 \markCellD{4}{6}{8}{9}
\end{scope}
\end{tikzpicture}
  \caption{A $2\times2$ mesh with $Q_1$ finite elements. Left: Local degrees of 
           freedom in each cell. Right: Global degrees of freedom after the 
           application of algorithm~\ref{alg:local_to_global_dof_map}.}
  \label{fig:local_to_global_dof_mesh}
 \end{center}
 \end{figure}
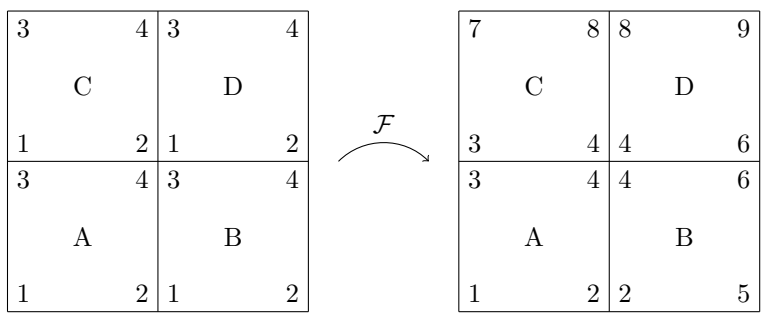
 \begin{figure}
  \newcommand{\drawGridEquivalenceRelation} %
{
 \draw[step=1.0,black!25,thin] (0,0) grid (16,1);
 \draw[step=4.0,black!75,thick] (0,0) grid (16,1);
 \draw[step=4.0,black!75,thick] (0,1) -- (16,1);

 \draw (1, 1) node[above right] {$A$};
 \draw (5, 1) node[above right] {$B$};
 \draw (9, 1) node[above right] {$C$};
 \draw (13,1) node[above right] {$D$};
}
\newcommand{\fillSquare}[2]{ \fill[#1] (#2,0) rectangle (#2+1,1); }
\newcommand{\putNumber}[1]{ \draw (-1.5,0.5) node {#1}; }

\begin{tikzpicture}[scale=0.3]
 \begin{scope}[xshift=300]
  \drawGridEquivalenceRelation
 \end{scope}

 \begin{scope}[yshift=-100]
  \fillSquare{fill=red!50, postaction={pattern=north east lines}}{1}
  \fillSquare{fill=red!50, postaction={pattern=north east lines}}{4}
  \fillSquare{fill=green!50, postaction={pattern=north west lines}}{3}
  \fillSquare{fill=green!50, postaction={pattern=north west lines}}{6}
  \drawGridEquivalenceRelation
  \putNumber{(i)};
 \end{scope}

 \begin{scope}[xshift=600, yshift=-100]
  \fillSquare{fill=red!50, postaction={pattern=north east lines}}{1}
  \fillSquare{fill=red!50, postaction={pattern=north east lines}}{4}
  \fillSquare{fill=blue!50, postaction={pattern=vertical lines}}{2}
  \fillSquare{fill=blue!50, postaction={pattern=vertical lines}}{8}
  
  \fillSquare{fill=green!50, postaction={pattern=north west lines}}{3}
  \fillSquare{fill=green!50, postaction={pattern=north west lines}}{6}
  \fillSquare{fill=green!50, postaction={pattern=north west lines}}{9}
  \drawGridEquivalenceRelation
  \putNumber{(ii)};
 \end{scope}
  
 \begin{scope}[yshift=-200]
  \fillSquare{fill=red!50, postaction={pattern=north east lines}}{1}
  \fillSquare{fill=red!50, postaction={pattern=north east lines}}{4}
  \fillSquare{fill=blue!50, postaction={pattern=vertical lines}}{2}
  \fillSquare{fill=blue!50, postaction={pattern=vertical lines}}{8}
  
  \fillSquare{fill=green!50, postaction={pattern=north west lines}}{3}
  \fillSquare{fill=green!50, postaction={pattern=north west lines}}{6}
  \fillSquare{fill=green!50, postaction={pattern=north west lines}}{9}
  \fillSquare{fill=green!50, postaction={pattern=north west lines}}{12}
  
  \fillSquare{fill=yellow, postaction={pattern=horizontal lines}}{7}
  \fillSquare{fill=yellow, postaction={pattern=horizontal lines}}{13}
  \drawGridEquivalenceRelation
  \putNumber{(iii)};
 \end{scope}

 \begin{scope}[xshift=600, yshift=-200]
  \fillSquare{fill=red!50, postaction={pattern=north east lines}}{1}
  \fillSquare{fill=red!50, postaction={pattern=north east lines}}{4}
  \fillSquare{fill=blue!50, postaction={pattern=vertical lines}}{2}
  \fillSquare{fill=blue!50, postaction={pattern=vertical lines}}{8}
  
  \fillSquare{fill=green!50, postaction={pattern=north west lines}}{3}
  \fillSquare{fill=green!50, postaction={pattern=north west lines}}{6}
  \fillSquare{fill=green!50, postaction={pattern=north west lines}}{9}
  \fillSquare{fill=green!50, postaction={pattern=north west lines}}{12}
  
  \fillSquare{fill=yellow, postaction={pattern=horizontal lines}}{7}
  \fillSquare{fill=yellow, postaction={pattern=horizontal lines}}{13}
  
  \fillSquare{fill=green!50!blue!50, postaction={pattern=dots}}{11}
  \fillSquare{fill=green!50!blue!50, postaction={pattern=dots}}{14}
  \drawGridEquivalenceRelation
  \putNumber{(iv)};
 \end{scope}
\end{tikzpicture}
  \caption{The partition $P$ of $M$ during the four relevant steps in 
           Algorithm~\ref{alg:local_to_global_dof_map}. Each cell corresponds
           to one local degree of freedom. Cells with the same color/pattern are
           considered equal, only empty (white) cells are not equal to each
           other.}
  \label{fig:local_to_global_dof_equivalence_relation}
 \end{figure}
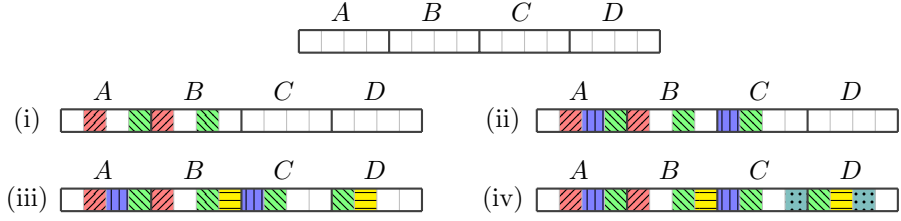
 The set of all local degrees of freedom is therefore 
 $M = \left\{ (A,1),\dots,(A,4) \right\} 
 \cup \left\{ (B,1),\dots,(B,4) \right\}
 \cup \left\{ (C,1),\dots,(C,4) \right\}
 \cup \left\{ (D,1),\dots,(D,4) \right\}$.
 Algorithm~\ref{alg:local_to_global_dof_map} modifies a partition $P$, 
 which is initialized to be the finest partition of $M$. During the algorithm,
 the condition $id(K) < id(K')$ is true exactly four times and steps 
 \ref{alg_line:find_partitioning} and \ref{alg_line:update_partitioning} are 
 depicted in Figure~\ref{fig:local_to_global_dof_equivalence_relation}.
 \begin{enumerate}[(i)]
   \item $K = A$, $K' = B$: Identify $(A,2),(B,1)$ (red) and $(A,4),(B,3)$
         (green).
   \item $K = A$, $K' = C$: Identify $(A,3),(C,1)$ (blue) and $(A,4),(C,2)$
         (green).
   \item $K = B$, $K' = D$: Identify $(B,3),(D,1)$ (green) and $(B,4),(D,2)$
         (yellow).
   \item $K = C$, $K' = D$: Identity $(C,2),(D,1)$ (green) and $(C,4),(D,3)$
         (teal).
  \end{enumerate}
 Now the set $P$ defines a partition with nine sets, each describing one global
 degree of freedom. The map $\mathcal F$ identifies $M$ with its image space
 $\{1,\dots,9\}$:
 \begin{align*}
   \mathcal F(A,1) & = 1, &&&
   \mathcal F(A,2) & = \mathcal F(B,1) = 2, \\
   \mathcal F(A,3) & = \mathcal F(C,1) = 3, &&&
   \mathcal F(A,4) & = \mathcal F(B,3) = \mathcal F(C,2) 
                       = \mathcal F(D,1) = 4, \\
   \mathcal F(B,2) & = 5, &&&
   \mathcal F(B,4) & = \mathcal F(D,2) = 6, \\
   \mathcal F(C,3) & = 7, &&&
   \mathcal F(C,4) & = \mathcal F(D,3) = 8, \\
   \mathcal F(D,4) & = 9.
  \end{align*}
 The image of $\mathcal F$ is illustrated in Figure~\ref{fig:local_to_global_dof_mesh} on the right.
\end{example}

The implementation of the concepts described in this section have been
adapted from {\sc MooNMD} to {\sc ParMooN} in a 
straightforward way. It is clear that the definition of the local basis and functionals  
is completely local and therefore not
affected by the parallelization. The required continuity in the definition of finite 
element spaces 
is performed in {\sc ParMooN} on each process separately.
The computed results are visualized with the software package {\sc ParaView} \cite{paraview_book,paraview_guide} that  does not 
require global numbers for the d.o.f.s. 
However, if needed, global numbers for d.o.f.s across all processes can be assigned, compare \cite{GJM+16}.
Altogether, the used concepts 
turned out to be  applicable in the same way  for the sequential as well as for the parallel code.

\section{Parallel Data Structures in {\sc ParMooN}}\label{sec:par_data}

In this chapter, the main steps in our approach to parallelize a finite element code are described.
{\sc ParMooN} supports a {single program, multiple data} (SPMD) approach on
parallelism using the Message Passing Interface (MPI) standard 
\cite{mpistandard2015}. It relies on a 
decomposition of the domain, which is the standard for parallelized finite 
element codes. Decomposing the computational domain and distributing it among 
the processes naturally leads to a parallelization of matrix-vector operations. 
The local aspect of the finite element method, which is reflected in the 
sparsity of the arising matrices, limits the communication overhead.

The main steps to be described in this section are domain decomposition, d.o.f. classification, consistency levels for distributed vectors, and technical details on the implementation of communication.

\subsection{Decomposing the Domain - Own Cells and Halo 
Cells}\label{sec:decomp}

In order to distribute the domain among the participating MPI 
processes,  {\sc ParMooN} makes use of the METIS graph partitioning tool \cite{METIS}. 
At program start, 
all processes read the same geometry and perform the same 
initial domain refinement steps. Upon reaching the first refinement level on 
which  computations will be performed, the root process
(process number 0) calls the {\sc METIS} library to compute a disjoint domain 
decomposition, i.e., to determine which process is going to be in charge of which 
mesh cells.

\begin{figure}[t!]
\centering
\includegraphics[width=0.9\textwidth]{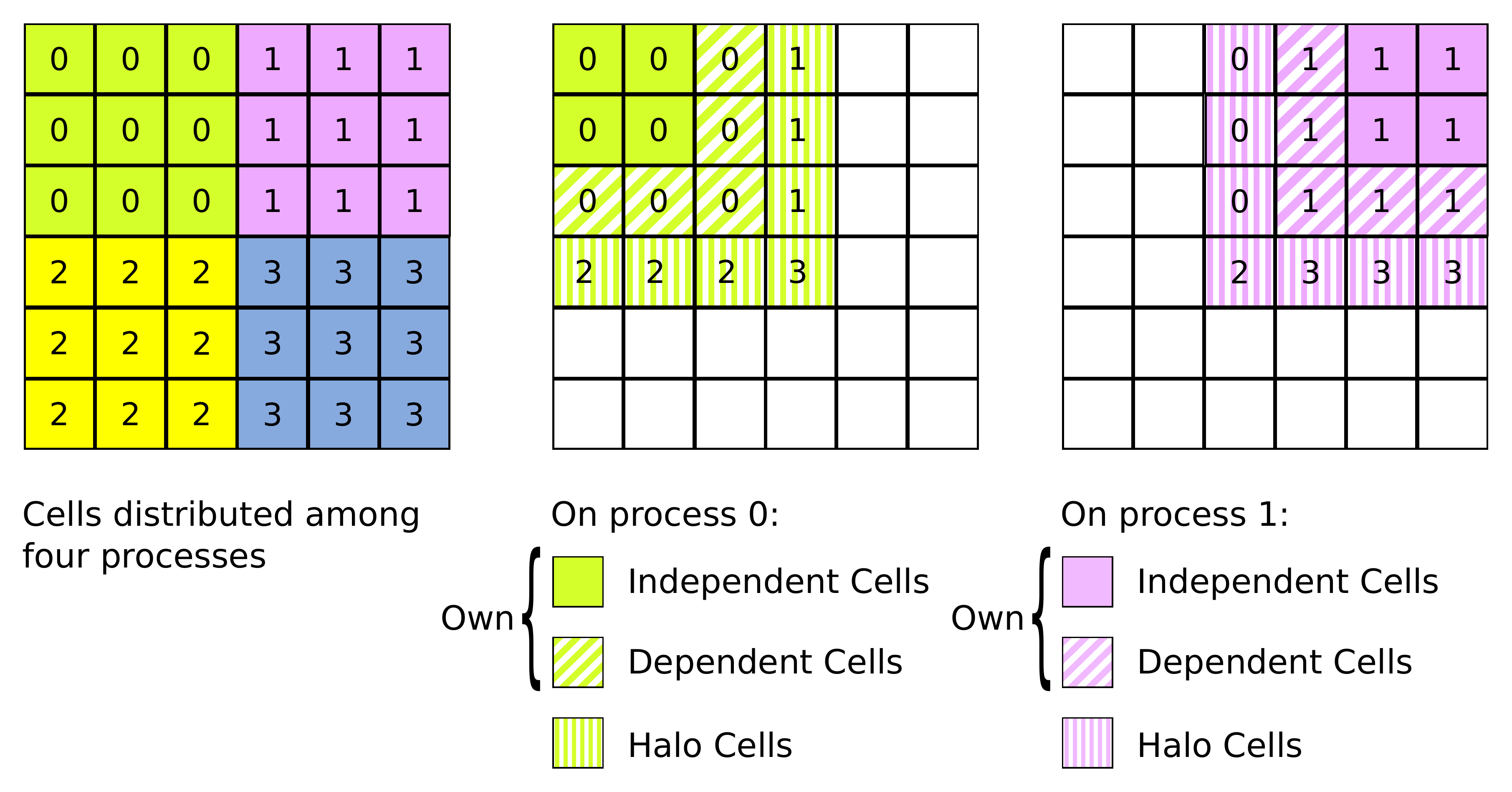}
\caption{Different cell types due to a domain decomposition.}
\label{cells_fig}
\end{figure}

Root then communicates the {\sc METIS} output to the other processes. Each process is 
informed about the cells it will be responsible for. These cells are 
called \textit{own cells} of the process. 
Each process then keeps only its own cells plus those cells that share a boundary 
face, edge, or vertex with an own cell. In domain decomposition methods these 
cells are commonly referred to as \textit{halo cells}. The sketches in Figure~\ref{cells_fig} 
clarify that expression -- the halo 
cells form a one-layer thick halo around the set of own cells. Each process $P$ deletes all cells which are neither own nor halo cells of $P$.

The own cells are further divided into \textit{dependent} and 
\textit{independent cells}. Therefore one defines the \textit{interface} as the set of those faces, edges and vertices which are shared by own cells and halo cells. All
own cells that contain a piece of interface are called dependent cells, while the remaining own cells 
are called independent cells.

The requirements on an efficient domain decomposition are twofold: the computational 
load must be balanced, i.e., there should be a comparable number of cells on each process, and the 
needed amount of communication must be small, i.e., the interfaces should be as small as possible.
Due to the deletion of cells, each process stores only a part of the entire
problem, but all processes execute the same program code. This was initially referred to as the single program, multiple data approach. 

With its domain reduced to own cells plus halo cells, each process sets up a 
finite element space and performs
all further computations only on its part of the domain.

\subsection{Types of Degrees of Freedom}

In the case of parallelized finite element methods, communication is necessary to 
transmit values of d.o.f.s across processes. 
This section describes the classification of d.o.f.s which is applied for this purpose.

A d.o.f. 
with index $i$ is defined by the finite element basis function 
$\phi_i$ and the associated global linear functional $\Phi_i$. It is represented, e.g., by the $i$-th entry
in the vectors for the solution and the right-hand side.
Each process is responsible for the classification of its d.o.f.s. The class to which a certain d.o.f. $i$ belongs to depends on its location and on the classes of those d.o.f.s with which it is coupled.
That is to say, two degrees of freedom of a finite element in $d$ dimensions are said to be \textit{coupled} if
the supports of the corresponding basis functions intersect on a set of non-zero $d$-dimensional 
measure. Note that, with this definition, coupling only occurs for d.o.f.s that are located in the 
same mesh cell. This notion of coupling transfers directly to a property of the finite element matrix $A$: the 
coupling of d.o.f. 
$i$ and $j$ will (potentially) lead to non-zero entries $a_{ij}$ and $a_{ji}$, and the needed memory 
for these entries has to be allocated.

In the following, the classes handling d.o.f.s in {\sc ParMooN} will be described shortly. 
All d.o.f.s that are localized in a own cell known to process $P$ are called \textit{known} d.o.f.s and their set will be denoted by $D^P_{\mathrm{known}}$ ($D^P_{\ast}$  stands  in the following for the set of d.o.f.s on $P$ of type $\ast$).
This set is then divided into \textit{masters} and \textit{slaves}, i.e., 
\[
D^P_{\text{known}} = D^P_{\text{master}} \overset{\Bigcdot}{\cup} D^P_{\text{slave}},
\]
where $\overset{\Bigcdot}{\cup}$ denotes the disjoint union.
A d.o.f. $i$ is said to be a master on $P$ if  $P$ is responsible 
for the value of $i$, in a way that will be clarified in Section \ref{consistency_levels}.
All known d.o.f.s which are not 
master on $P$ are called slaves. It is worth noting that every d.o.f. in the 
entire problem is  a master on exactly one process.  In contrast, a d.o.f. can be slave on more 
than one process. 

An even finer classification of $D^P_{\text{master}}$ and $D^P_{\text{slave}}$ will be used in the following.
Note
that the names for the classes correspond only loosely to the types of cells they are located in.
Care must be taken of those d.o.f.s that lie on the intersection of different cell types. The classes are:

\begin{list}{}{\itemsep0.0ex\parsep0.1ex\topsep0.0ex}
\item[$\bullet$] \textit{Independent d.o.f.}, i.e., 
all d.o.f.s which lie in $P$'s own cells but not in its dependent cells.
All $P$'s independent d.o.f.s are set as 
masters, since they are not even known to any other process. They only 
couple to other masters of $P$.
\item[$\bullet$]
\textit{Dependent d.o.f.}, i.e., 
those d.o.f.s lying in $P$'s dependent cells, but not in its halo cells. Process $P$ is the master
of all its dependent d.o.f.s.
The notation is motivated by the fact that {the dependent d.o.f.s} 
are in a vicinity to the domain interface and therefore possess a certain dependency 
on other processes. 
\item[$\bullet$]
\textit{Interface d.o.f.}, i.e., 
all d.o.f.s that lie on the intersection of dependent cells and halo cells.
These d.o.f.s are known to all adjacent 
processes as interface d.o.f.s, too. Only one of these processes will take master 
responsibility for each interface d.o.f. In particular, on a process $P$, one distinguishes between
\textit{master interface d.o.f.s} (the interface 
d.o.f.s which are master on $P$) and \textit{slave interface d.o.f.s} (all others d.o.f.s,  
for each of which a neighboring process takes master responsibility).
\item[$\bullet$] \textit{Halo d.o.f.}, i.e.,
all d.o.f.s which lie in halo cells but not on the interface.
Since all of them are dependent d.o.f.s to neighboring 
processes, one of these will take master responsibility for them. On $P$ all 
halo d.o.f.s are slaves.
\end{list}
Hence, the sets of masters and slaves can be divided as
\[
\begin{aligned}
& D^P_{\text{master}}  = D^P_{\text{independent}} \overset{\Bigcdot}{\cup} 
D^P_{\text{dependent}} 
\overset{\Bigcdot}{\cup} D^P_{\text{interface master}}\,,\\
& D^P_{\text{slave}} = D^P_{\text{interface slave}} \overset{\Bigcdot}{\cup} 
D^P_{\text{halo}}.
\end{aligned}
\]

\begin{figure}[t!]
\centering
\includegraphics[width=0.8\textwidth]{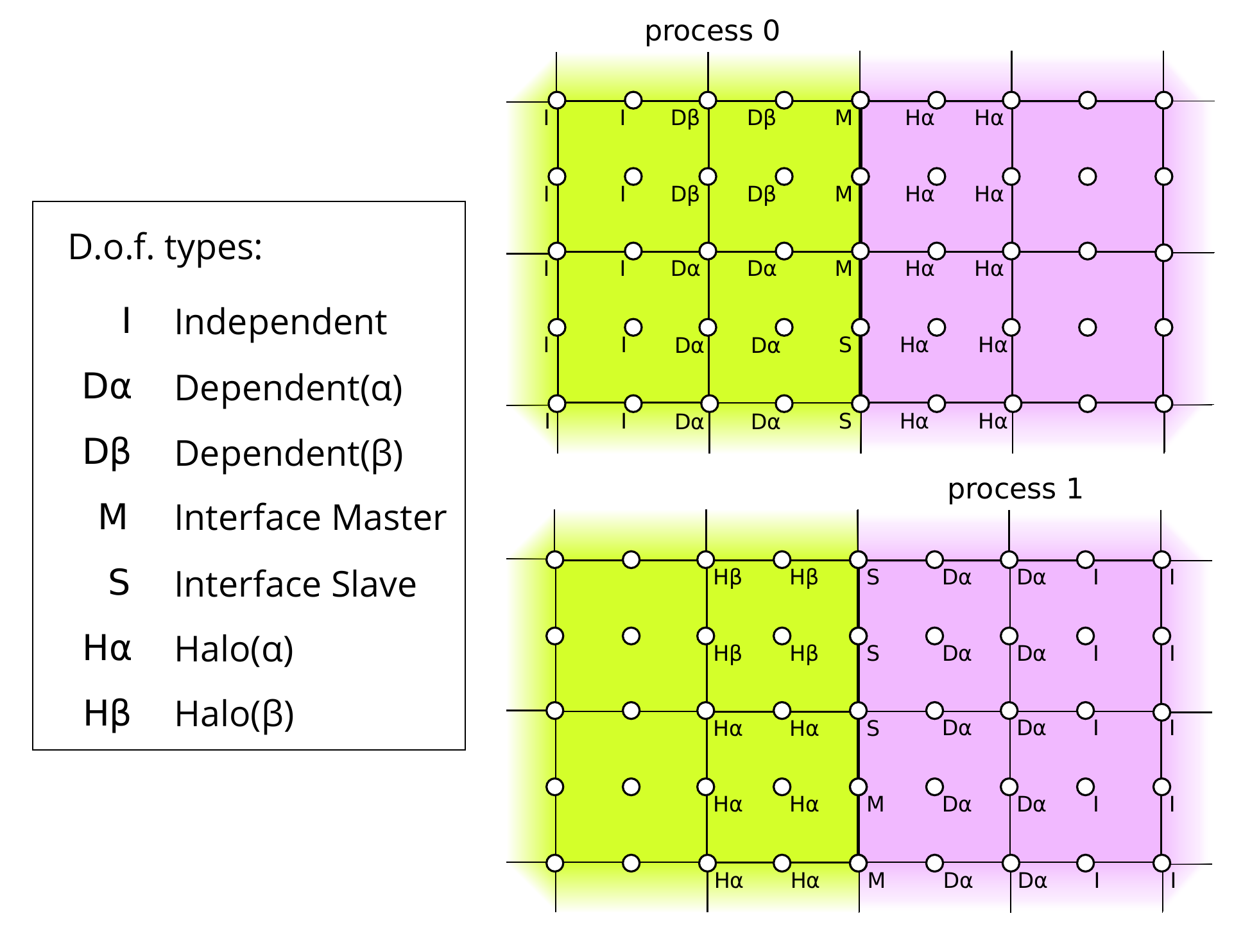}
\caption{Types of d.o.f.s at the interface for two-dimensional $Q_2$ finite 
elements, from the point of view of process~0 and process~1, where only the named d.o.f.s are known to the respective process.}
\label{dof_types}
\end{figure}

Moreover, it is convenient to further refine the d.o.f. classification in 
order to reduce the communication overhead, see Section~\ref{consistency_levels}.
To this aim, the halo and the dependent d.o.f.s of $P$ are further 
divided into Halo($\alpha$) and Halo($\beta$), and Dependent($\alpha$) and Dependent($\beta$), respectively:
\begin{list}{}{\itemsep0.0ex\parsep0.1ex\topsep0.0ex}
\item[$\bullet$]
\textit{Halo($\alpha$)} and \textit{Halo($\beta$) d.o.f.s.}:
The  Halo($\alpha$) d.o.f.s are 
those that are coupled with at least one (interface) master of $P$, while 
Halo($\beta$) d.o.f.s are coupled solely with other slaves, i.e., with interface slave d.o.f.s 
and other halo d.o.f.s.
\item[$\bullet$]
\textit{Dependent($\alpha$)} and \textit{Dependent($\beta$) d.o.f.s.}:
Dependent($\alpha$) d.o.f.s are those
connected to at least one (interface) slave, while Dependent($\beta$) d.o.f.s
are all those that are connected to masters only, i.e., to  
interface master, other dependent, or independent d.o.f.s.
Note that all Dependent($\beta$) d.o.f.s of process $P$ will be Halo($\beta$) on 
all other processes where they are known. For Dependent($\alpha$) d.o.f.s the situation is 
not as simple. Each of them is Halo($\alpha$) to at least one neighboring process, but 
can be Halo($\beta$) to others.
\end{list}

In general, the d.o.f. classification depends on the finite element spaces and the decomposition
of the domain. 
Figure~\ref{dof_types} sketches the different types for two-dimensional $Q_2$ finite 
elements at the interface between two processes.

\subsection{Consistency Levels}\label{consistency_levels}

In parallel computations, values and information, e.g., finite element vectors, can be stored either in a consistent way or 
in an additive way. In the case of a finite element vector, which is of major interest here, \textit{consistent} storage means 
that all processes have the same and correct value at all respective d.o.f.s 
as in the sequential environment. 
If the vector is present in \textit{additive} storage, each global value is the 
sum of the values over all processes where it is known.  In {\sc ParMooN}, certain concepts of weakened 
consistent storage of finite element vectors and matrices are used. 

A finite element vector belongs to one of the following four levels of consistency.
\begin{list}{}{\itemsep0.0ex\parsep0.1ex\topsep0.0ex}
\item[$\bullet$] \textit{Level-0-consistent.}  Consistency holds only with regard to 
       master degrees of freedom. Each master on each process holds the same value as 
       it would be in a sequential computation. The values of slaves are in 
       an undefined storage state. In the implementation of operations care must be taken not to lose
       Level-0-consistency -- all master values must be kept as they would be in a sequential execution.
 \item[$\bullet$]  \textit{Level-1-consistent.} All masters 
       and all interface slave d.o.f.s are stored consistently. The values of all halo 
       d.o.f.s are in an undefined storage state.
 \item[$\bullet$]  \textit{Level-2-consistent.}  Consistency is established for all but 
       Halo($\beta$) d.o.f.s. The values of Halo($\beta$) d.o.f.s. are in an undefined storage state
       while all other values are consistent.
 \item[$\bullet$]  \textit{Level-3-consistent.}  All  d.o.f.s are stored consistently. 
       This situation is called the ``real'' consistent storage.
\end{list}
The main motivations behind introducing this classification are that 
different operations require a different state of consistency of their input 
data  and that restoring a certain state of consistency requires a certain amount of 
communication -- the lower the required state of consistency, the lower the 
required amount of communication.

After the domain has been decomposed, each process $P$ assembles a finite element matrix 
on all its known cells. The use of halo cells assures that all 
information to assemble the rows belonging to masters is available on 
$P$.  The  finite element matrix assembled in this way will have the property that all rows and columns 
which belong to masters
are correctly assembled, i.e., consistent. A matrix with this property is called Level-0-consistent, too. In fact, even the rows belonging to interface slave d.o.f.s are correct, but this is only a by-product. Hence,
the finite element matrix is even Level-1-consistent. Note that it is not possible to extend this concept to Level-2- or even Level-3-consistency of a matrix. The reason 
is that some of the d.o.f.s which would be needed to store the entire matrix row associated with a halo d.o.f. are not in $D^P_{\text{known}}$.

Regarding operations, multiplication of a Level-0-consistent matrix with a 
Level-2-con\-sis\-tent vector gives a Level-0-consistent vector.  
Level-3-consistency of the input vector is not needed. If the input is 
Level-3-consistent, together with the matrix being Level-1-consistent, 
the result will be Level-1-consistent.

Multiplying a vector with a constant scalar will maintain the current 
consistency level, as will vector-vector addition. In the latter case, vectors 
with different consistency level might be added and the result has the lower 
consistency level. 
Scalar products require 
Level-0-consistency of both vectors, where all slaves will be skipped,  
and a globally additive reduce operation is required to get a consistent result.

Level-3-consistency of a finite element vector in {\sc ParMooN} is only enforced
if operations require knowledge of the represented finite 
element function even on the 
halo cells. Such operations include 
the matrix assembling with 
an input finite element function, e.g., for the convective term of the Navier--Stokes equations,
interpolating initial conditions in a time-dependent problem,  
 or gradient recovery by averaging 
gradients over a patch of mesh cells.

Enforcing certain consistency levels is a matter of communication. For each 
d.o.f. which needs an update, the responsible master process 
communicates its value to all processes where it is slave.  These processes simply 
reset its value to the received value. The required infrastructure for this communication 
is set up just once for a certain finite element space and can be reused whenever an 
update is necessary. 

\subsection{Organizing Communication}


When setting up the communication structure, one has to find, for each non-independent master $i$,  
all those slaves on other processes that are 
globally identical to $i$. Certain master types match with certain slave types, see 
Table~\ref{ms_relation}, forming three distinct pairs of master-slave 
relations. To restore a certain consistency level, at least one of these relations requires an 
update.

\begin{table}[t!]
\begin{center}
\caption{Master-slave relationship of d.o.f. types.}
\begin{tabular}{l l c l}
\textbf{Relation (shorthand)} & \textbf{Master type} & updates &\textbf{Slave type} \\ 
\hline
Interface (IMS) & Interface Master & $\longrightarrow$ & Interface Slave \\ 
Dependent($\alpha$)--Halo($\alpha$) (DH$\alpha$) & Dependent($\alpha$) &$\longrightarrow$ & Halo($\alpha$) \\ 
Dependent($\beta$)--Halo($\beta$) (DH$\beta$) & Dependent($\beta$) & $\longrightarrow$& Halo($\beta$) \\ 
\end{tabular}
\label{ms_relation}
\end{center}
\end{table}

Note that it is generally not possible to immediately identify the global number of a
d.o.f., since each process creates finite element spaces only on its 
known (own and halo) cells and it numbers its d.o.f.s locally, unaware of the other processes.
In {\sc ParMooN}, the global identification of a d.o.f. is defined according to 
the global 
number of the mesh cell in which the d.o.f. is located.
These global cell numbers are assigned to 
each cell before decomposing the domain. Performing only uniform refinements after domain decomposition, 
such a globally unique cell number can easily be given to 
children cells, too. 
Hence, the global cell number 
and a consistent local numbering of the d.o.f.s within each cell enable to 
identify each d.o.f. globally.

The communication structure is stored in a class called {\acrfont ParFEMapper}, while the 
communication itself is performed by a class named {\acrfont ParFECommunicator}. 
Setting up the {\acrfont ParFEMapper} and {\acrfont ParFECommunicator} requires some 
communication itself, and a detailed description of this task is presented in \cite{GJM+16}. 
In what follows, only an overview of it is given, including a short description of those data fields of 
{\acrfont ParFEMapper} that are relevant when updating the d.o.f.s of a certain 
master-slave relation, see Table~\ref{ms_relation}.
The interface (IMS) relation is presented as an example. For the two other relations, the data fields are defined in a similar way.

For the IMS update, the {\acrfont ParFECommunicator} wraps a call to the MPI function 
{\acrfont MPI\_Alltoallv}, where every process may send a different 
set of values of the same type ({\acrfont MPI\_DOUBLE} in this case) to each 
other process. To control the {\acrfont MPI\_Alltoallv} call, the {\acrfont ParFEMapper}
stores the following data, where {\acrfont mpi\_size} is the total number of 
processes and {\acrfont nInterfaceSlaves} is the number of interface slaves 
{local} to process $P$.
\begin{list}{}{\itemsep0.0ex\parsep0.1ex\topsep0.0ex}
 \item[$\bullet$] {\acrfont int* sendBufIMS}  is the send buffer, filled with the values 
of all interface masters, each one possibly appearing more than once, 
which will then be sent to the other processes. Its total length equals the sum 
over all values of {\acrfont sendCountsIMS}.
 \item[$\bullet$] {\acrfont int* sendCountsIMS} is an array of size 
{\acrfont mpi\_size}. It lists how many values $P$ has to send to each 
process.
 \item[$\bullet$] {\acrfont int* sendDisplIMS} is the send displacement, an array of size 
{\acrfont mpi\_size}. It stores where in the array {\acrfont sendBufIMS} 
the message for a certain process starts. In {\sc ParMooN} there are neither 
overlaps nor gaps, so {\acrfont sendDisplIMS[i]} holds 
the sum of {\acrfont sendCountsIMS[0]} to {\acrfont sendCountsIMS[i-1]}.
 \item[$\bullet$] {\acrfont int* recvBufIMS} is the receive buffer, which will be filled with 
sent values from the other processes in the communication routine. Its size 
equals {\acrfont nInterfaceSlaves}.
 \item[$\bullet$] {\acrfont int* recvCountsIMS} is an array of  size {\acrfont 
mpi\_size}. It lists how many values are to be received from each 
process. The sum of its values is {\acrfont nInterfaceSlaves}.
 \item[$\bullet$] {\acrfont int* recvDisplIMS} is  the receive displacement. Like for 
 {\acrfont int* sendDisplIMS} there are neither gaps nor overlaps. 
\end{list}
Besides the data, which is needed in the immediate control of {\acrfont 
MPI\_Alltoallv}, the {\acrfont ParFEMapper} contains two arrays that allow to interpret the 
sent and received data, by mapping between send- or receive buffer and the local 
d.o.f.s:
\begin{list}{}{\itemsep0.0ex\parsep0.1ex\topsep0.0ex}
 \item[$\bullet$]  {\acrfont int* sentDofIMS} interprets {\acrfont sentDofIMS[i]}$= d$ as: 
The $i$-th place in {\acrfont sendBufIMS} has to be filled 
with the value of the local d.o.f. $i$.
 \item[$\bullet$]  {\acrfont int* rcvdDofIMS} interprets {\acrfont rcvdDofIMS[i]}$= d$ as: 
The $i$-th value in {\acrfont recvBufIMS} should update the
local d.o.f. $i$.
\end{list}

To change a vector from Level-0-consistency to Level-1-consistency, 
only an
{IMS} update is required. For restoring Level-2-consistency, 
additionally a {DH$\alpha$} update is necessary, while Level-3-consistency 
requires even a {DH$\beta$} update on top. 

\section{The Parallel Geometric Multigrid Method}\label{sec:par_mg}

Geometric multigrid methods are an appealing option to be used as preconditioners 
in problems where the necessary hierarchy of grids can be provided. These methods 
have been used in the simulations performed with {\sc MooNMD} in particular for 
three-dimensional problems and for linear saddle point problems arising in the 
linearization and discretization of equations modeling incompressible flow 
problems \cite{Joh16}.

The components of a geometric multigrid method are the following: function prolongation, 
defect restriction, function restriction, smoother, and coarse grid solver. Concerning details 
of the algorithms and implementation of the first three components in {\sc MooNMD}/{\sc ParMooN}, 
it will be referred to \cite{JM04} for details. 
It shall be only noted that multilevel methods are 
supported that allow different finite element spaces on different levels of the 
multigrid hierarchy. In particular, the so-called multiple discretization multilevel
(MDML) method for higher order discretizations can be used. 
This method possesses  the higher order discretization on the highest multigrid level
and a low order discretization on all coarser levels, where  more multigrid levels 
might be defined than geometric levels exist. The motivation for this 
approach is the experience that multigrid methods often work very efficiently for low
order discretizations. Numerical studies of the efficiency of the MDML method can 
be found, e.g., in \cite{Joh02}. The grid transfer operations are performed 
with a local operator, taking values only on the mesh cells of the current level, 
proposed in \cite{Sch00}, that 
can handle different finite element spaces on different levels of the 
grid hierarchy, see \cite{Joh16,JM04} for details.

The implementation and parallelization of geometric multigrid methods require a 
considerable amount of work. The geometric data structures need to be equipped with 
parent-child information and the grid transfer operations have to be implemented. 
In the current version of {\sc ParMooN}, each process is responsible for a part of the coarsest 
grid, compare Section~\ref{sec:decomp}, and refines this part uniformly. Consequently, 
all parent-child information is available on the process.
More technical details on constructing the grid hierarchy can be found in \cite{GJM+16}.

As usual in parallel geometric multigrid 
methods, block-Jacobi smoothers are applied, where the blocks correspond to the 
master and interface slave  degrees of freedom of a process. Within the blocks, the actual 
smoother, like SSOR or the Vanka smoother, is used. After each smoothing iteration, the 
values at the interface are updated by computing their arithmetic average. 
As already mentioned above, 
the grid transfer operators need as potential input all values that are connected 
to a mesh cell. Since it is sufficient to perform the grid transfer only on own cells, 
the input vectors for the grid transfer have to be Level-1-consistent. 

\section{Numerical Studies}\label{sec:numres}

The performed numerical studies are a first step of assessing the efficiency of the parallelized geometric
multigrid method in comparison with parallel solvers that can be used by linking an external 
library to the code. The underlying question is whether it was worthwhile to perform 
the complex parallelization of this method. We think that this question arises also 
in other groups that maintain an in-house research code. 

In this paper, the numerical studies concentrate on the standard
multigrid method (same number of geometric and multigrid levels, same discretization
on each level) since we think that this method is of most interest for the community.
It was used as preconditioner in the flexible GMRES (FGMRES) method 
\cite{Saa93}.
The system on the 
coarsest grid was solved with a sparse direct solver. Moreover, 
the V-cycle was applied because of its better efficiency on parallel computers
compared with the more stable F- and W-cycle (the F-cycle is our standard approach in sequential simulations).
The V-cycle  approaches less often coarser grids, which possess an unfavorable ratio of 
computational work and necessary communications, than the other cycles.

The efficiency of the geometric
multigrid preconditioner has been compared with the efficiency of  
the sparse direct solver 
{\sc MUMPS} \cite{ADLK01,AGLP06}  and the  FGMRES method preconditioned 
with SSOR (for scalar problems), the BoomerAMG \cite{HY02} (for scalar problems), or the LSC preconditioner \cite{ESW14} (for linear saddle point problems). These 
solvers were used as they are provided in the library {\sc PETSc}, version 3.7.2, 
\cite{petsc-web-page,petsc-user-ref,petsc-efficient}. The restart parameter in FGMRES 
was set to be $50$.

The numerical studies were performed on a hardware platform that can be 
usually found in universities and academic institutes, in our case,  on compute servers HP BL460c Gen9
2xXeon, Fourteen-Core 2600MHz. 
We think that the performance 
on a hardware platform that is widely available is of interest for the community.
The results will consider only the computing times for the different solvers of the 
linear systems of equations. We could observe some  variations of these 
times for the same code and input parameters but different runs. To reduce the 
influence of these variations, all simulations were performed five times, 
the fastest and the slowest computing time were neglected 
and the average of the remaining three times is presented below. 

\begin{example} \rm {\em Steady-state convection-diffusion equation.}\label{ex:hemker}
This example is a three-dimensional extension of a benchmark problem for two-dimensional
convection-diffusion equations -- the so-called Hemker example \cite{ACF+11,Hem96}.
The domain is defined by 
$$
\Omega = \left\{\left\{(-3,9) \times (-3,3) \right\}\setminus\left\{ (x,y)\ : \ 
x^2+y^2 \le 1 \right\}\right\} \times (0,6)
$$
and the equation is given by 
$$
- \varepsilon \Delta u + \boldsymbol b \cdot \nabla u =  0 \quad \mbox{in } \Omega, 
$$
with $\varepsilon = 10^{-6}$ and $\boldsymbol b = (1,0,0)^T$.  Dirichlet boundary
conditions $u=0$ were prescribed at the inlet plane $\{x=-3\}$ and $u=1$ at the cylinder. 
At all other boundaries, homogeneous Neumann boundary conditions were imposed. 
This example models, e.g., the heat transport from the cylinder. The 
solution exhibits boundary layers at the cylinder and internal layers downwind the 
cylinder, see Figure~\ref{fig:hemker}.

\begin{figure}[t!]
\centerline{\includegraphics[height=0.23\textheight]{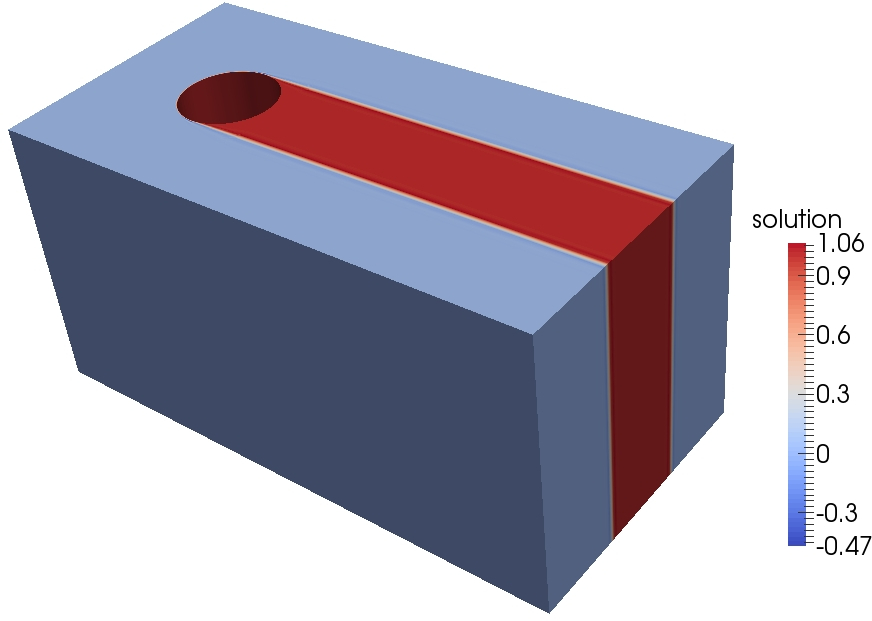}
\hspace*{0.5em} \includegraphics[height=0.23\textheight]{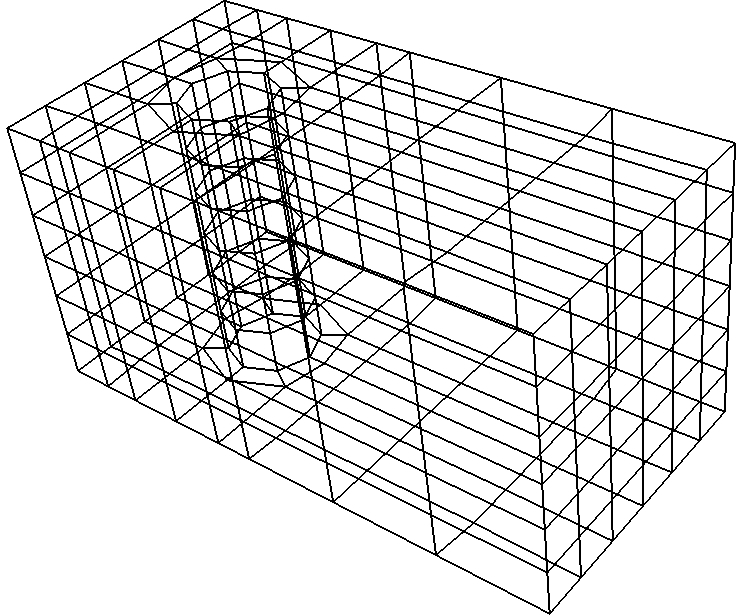}}
\caption{Example~\ref{ex:hemker}, numerical solution (left) and initial grid (level~0, right). The 
color bar shows that the numerical solution computed with the SUPG method possesses 
under- and overshoots. These spurious oscillations occur in particular in a vicinity 
of the cylinder, compare \cite[Fig.~14]{JKS11} for the two-dimensional situation.}
\label{fig:hemker}
\end{figure}

It is well known that stabilized discretizations have to be employed in the presence 
of dominant convection. In the numerical studies, the popular streamline-upwind
Petrov--Galerkin (SUPG) method \cite{BH82,HB79} was used with the standard parameter
choice given in \cite[Eqs.~(5) -- (7)]{JK07}. Simulations were performed 
with $Q_1$ finite elements.
The initial grid is depicted in Figure~\ref{fig:hemker}.

Computational results are presented in Figures~\ref{fig:hemker_results_4} 
and~\ref{fig:hemker_results_5}
for refinement levels~4 (1~297~375 d.o.f.s) and~5 (10~388~032 d.o.f.s). The {\sc PETSc} solvers were 
called with the flags {\tt -ksp\_type fgmres -pc\_type sor}
and {\tt -ksp\_type fgmres -pc\_type hypre -pc\_hypre\_type boomeramg}, respectively. 
For the geometric multigrid preconditioner, the V(2,2)-cycle was used and the overrelaxation parameter 
of the SSOR smoother within the block-Jacobi method was set to be $\omega=1$. 
This approach is certainly not optimal since with an increasing number of processes the number 
of blocks of the block-Jacobi method increases, which in turn makes it advantageous to apply 
some damping, i.e., a somewhat smaller overrelaxation parameter. In order not to increase the 
complexity of the numerical studies, we decided to fix a constant overrelaxation parameter 
that worked reasonably well for the whole range of processor numbers which was used. 
The iterative solvers were stopped if the Euclidean norm of the residual vector was 
less than $10^{-10}$. Parameters like the stopping criterion, the overrelaxation factor, and the restart parameter
were the same in all iterative methods.

\begin{figure}[t!]
\centerline{\includegraphics[width=0.48\textwidth]{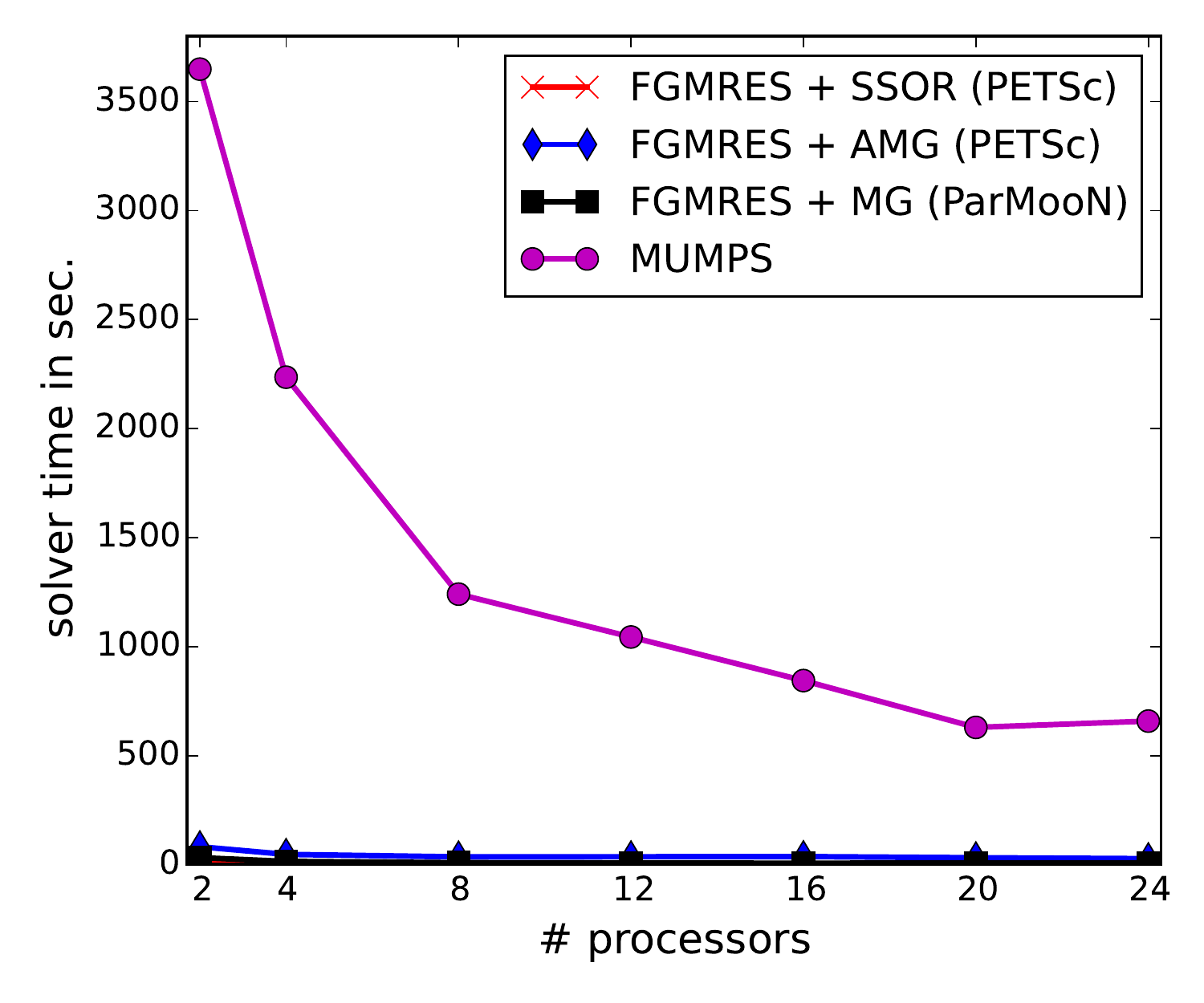}
\hspace*{1em} \includegraphics[width=0.48\textwidth]{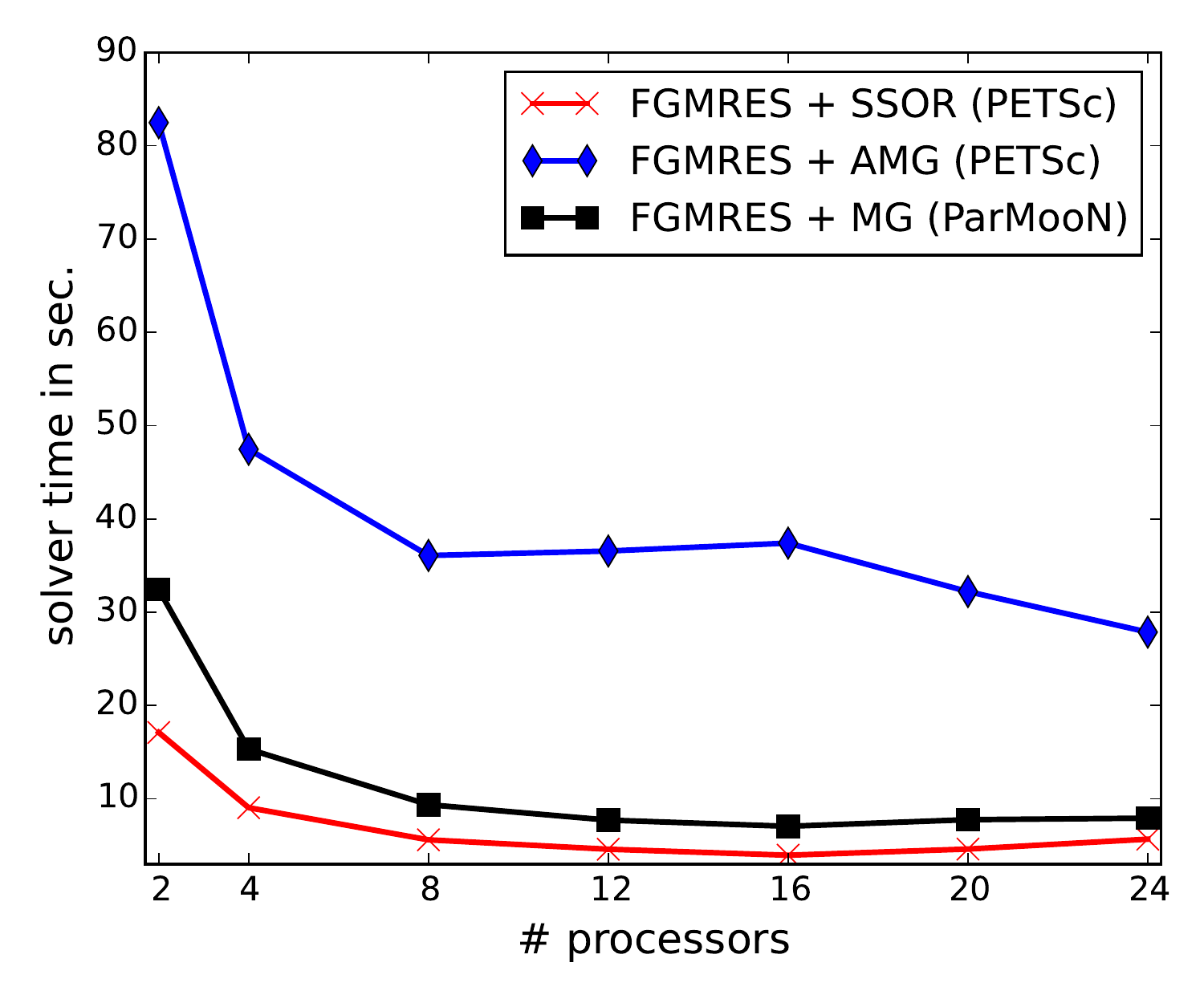}}
\caption{Example~\ref{ex:hemker}, solver times on refinement level~4.}
\label{fig:hemker_results_4}
\end{figure}
\begin{figure}[t!]
\centerline{\includegraphics[width=0.48\textwidth]{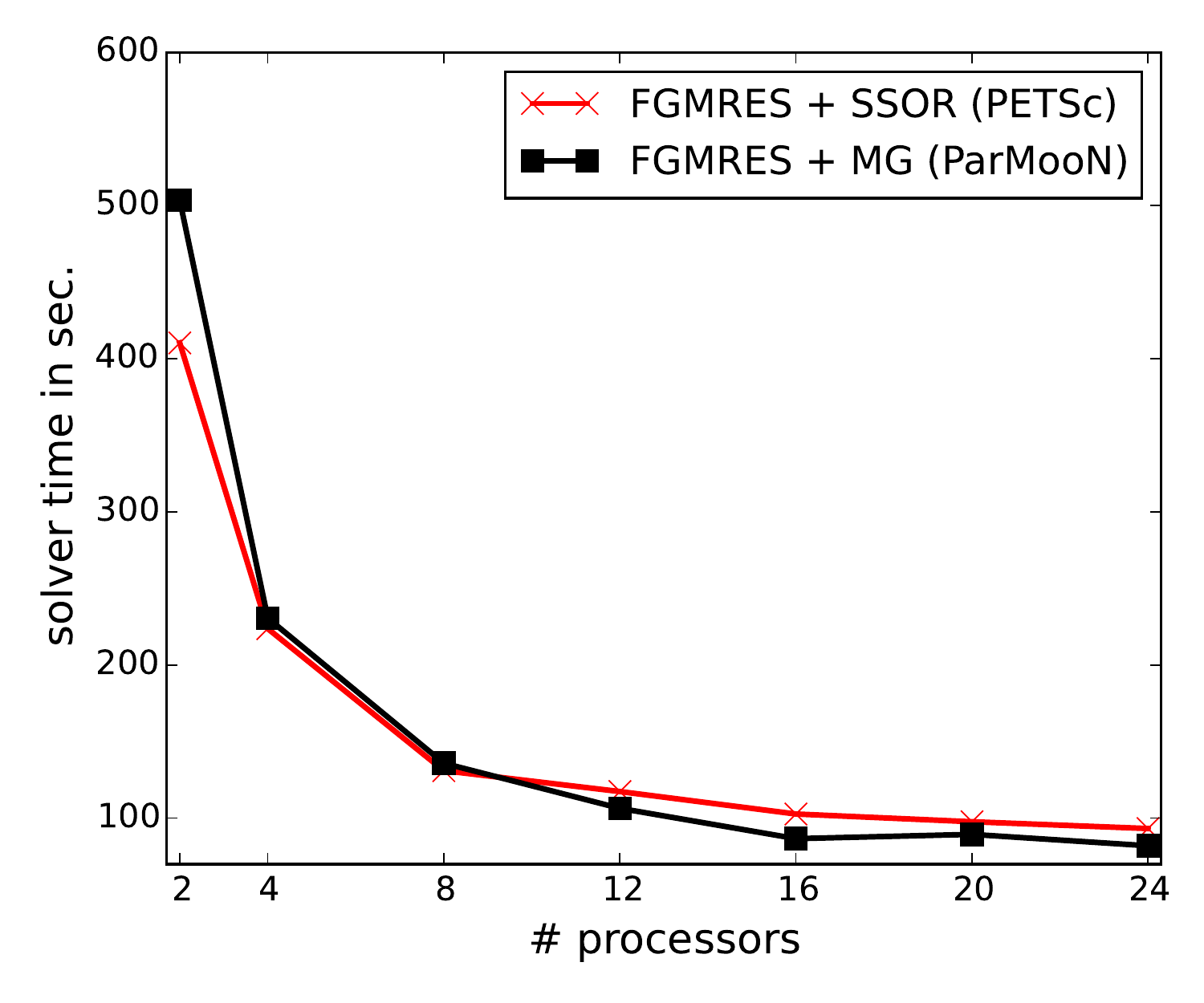}}
\caption{Example~\ref{ex:hemker}, solver times on refinement level~5. {\sc PETSc} FGMRES with  BoomerAMG converged only 
with two processors (1519~sec.).}
\label{fig:hemker_results_5}
\end{figure}

It can be seen in Figure~\ref{fig:hemker_results_4} that the sparse direct solver performed less efficient 
by around two orders of magnitude than the other solvers. This behavior was observed for all studied scalar problems
and no further results with this solver for scalar problems will be presented. Among the iterative solvers, 
{\sc PETSc} FGMRES with SSOR preconditioner and {\sc ParMooN} FGMRES with geometric multigrid preconditioner (MG) performed notably  more efficient
than {\sc PETSc} FGMRES with  BoomerAMG. The latter solver did not even converge on the finer grid in the simulations
on more than two processors. With respect to the computing times and the stability, similar results 
were obtained by applying two V(1,1) BoomerAMG cycles ({\tt -pc\_hypre\_boomeramg\_max\_iter 2}).
With the V(2,2) cycle ({\tt -pc\_hypre\_boomeramg\_grid\_sweeps\_all 2}), the solver did not converge. 
FGMRES with SSOR required considerably more iterations than {\sc ParMooN} FGMRES 
with MG: around 125 vs. 25 on level~4 and 320 vs. 45 on level~5. A notable decrease of the computing time for the iterative
solvers on the coarser grid can be observed only until 8~processors. On this grid, {\sc PETSc} FGMRES with SSOR was 
a little bit faster than {\sc ParMooN} FGMRES with MG. On the finer grid, a decrease of the computing times occurred 
until 16 processors and {\sc ParMooN} FGMRES with MG was often a little bit more efficient. 
\end{example}

\begin{example}\rm {\em Time-dependent convection-diffusion-reaction equation.}
\label{ex:bail}
This example can be found also in the literature, e.g., in \cite{JN12}, and 
it models a typical situation which is
encountered in applications. A species enters the domain $\Omega=(0,1)^3$ at the inlet $\Gamma_\text{in} = \{0\}\times(5/8,6/8)\times(5/8,6/8)$ and it is
transported through the domain to the outlet 
$\Gamma_\text{out} = \{1\}\times(3/8,4/8)\times(4/8,5/8)$. In addition, the
species is diffused somewhat and in the subregion where the species
is transported, also a reaction occurs. The ratio of
diffusion and convection is typical for many applications. 

The underlying model is given by 
\begin{equation*}
\begin{alignedat}{2}
\partial_t u-\varepsilon \Delta u+\boldsymbol{b}\cdot \nabla u+c u & = 0 &\quad& \mbox{ in } (0,3) \times \Omega,\\
u & =u_\text{in}  && \mbox{ in } (0,3) \times \Gamma_\text{in},\\
\varepsilon \frac{\partial u}{\partial \boldsymbol{n}} &= 0 && \mbox{ on } (0,3) \times \Gamma_{\text{N}},\\
u &= 0 &&  \mbox{ on } (0,3) \times \left(\partial\Omega \setminus \left(\Gamma_{\text{N}}
\cup \Gamma_\text{in}\right)\right), \\
u(0,\cdot) & =u_0 && \text{ in }\Omega.
\end{alignedat}
\end{equation*}
The diffusion parameter is given by $\varepsilon=10^{-6}$, 
the convection field is defined by ${\boldsymbol b}=(1,-1/4,-1/8)^T$, and the reaction by 
\begin{equation*}
c(\boldsymbol{x})= \begin{cases}
1 & \text{if } \text{dist}(\boldsymbol{x},g)\le 0.1,\\
0 & \text{else},
\end{cases}
\end{equation*}
where $g$ is the line through the center of the inlet and the center of the
outlet and $\text{dist}(\boldsymbol{x},g)$ denotes the Euclidean
distance of the point $\boldsymbol{x}$ to the line $g$. 
The boundary condition at the inlet is prescribed by 
\begin{equation*}
u_\text{in}= \begin{cases}
\sin(\pi t/2) & \text{ if } t\in [0,1], \\ 
1 & \text{ if } t\in (1,2],\\
\sin(\pi(t-1)/2) & \text{ if } t\in (2,3].
\end{cases}
\end{equation*}
The initial condition is set to be $u_0(\bx) = 0$. Initially, in 
the time interval $[0,1]$, the inflow increases and the injected
species is transported towards the outlet.  Then, there is a constant inflow 
in the time interval
$(1,2]$ and the species reaches the outlet. At $t=2$, 
there is almost a steady-state solution. 
Finally, the inflow decreases in the time interval $(2,3]$, compare \cite{JN12}.

 The SUPG stabilization of the $Q_1$ finite element method 
was used and 
as temporal discretization, 
the Crank--Nicolson scheme with the equidistant time step $\Delta t = 10^{-2}$ was 
applied. Simulations were performed on grids with $64^3$ and $128^3$ cubic mesh cells. 
The coarsest grid for the geometric multigrid method possessed $4^3$ cubic mesh cells. 
The SSOR method was applied with the overrelaxation parameter $\omega  = 1.25$. As initial 
guess for the iterative solvers, the solution from the previous discrete time was used. 
The availability of a good initial guess and the dominance of the system matrix by the mass 
matrix are the main differences to the steady-state 
case. 

Computing times for the iterative solvers are presented in Figures~\ref{fig:bail_results_64}
and \ref{fig:bail_results_128}. On the coarser grid, a speed-up can be observed only until 8~processors
and on the finer grid until 12 to 16 processors. With respect to the efficiency of the iterative 
solvers, the same observations can be made as in the steady-state Example~\ref{ex:hemker}.
Using the V(2,2) cycle or two V(1,1) cycles in the BoomerAMG gave the same order of computing times
as the presented times of the V(1,1) cycle. 
The 
same statement holds true for the number of iterations per time step, e.g., on the finer grid there were up to 45 for 
{\sc PETSc} FGMRES with SSOR and 
around 1-3 for {\sc PETSc} FGMRES with  BoomerAMG and {\sc ParMooN} FGMRES with MG.

\begin{figure}[t!]
\centerline{\includegraphics[width=0.48\textwidth]{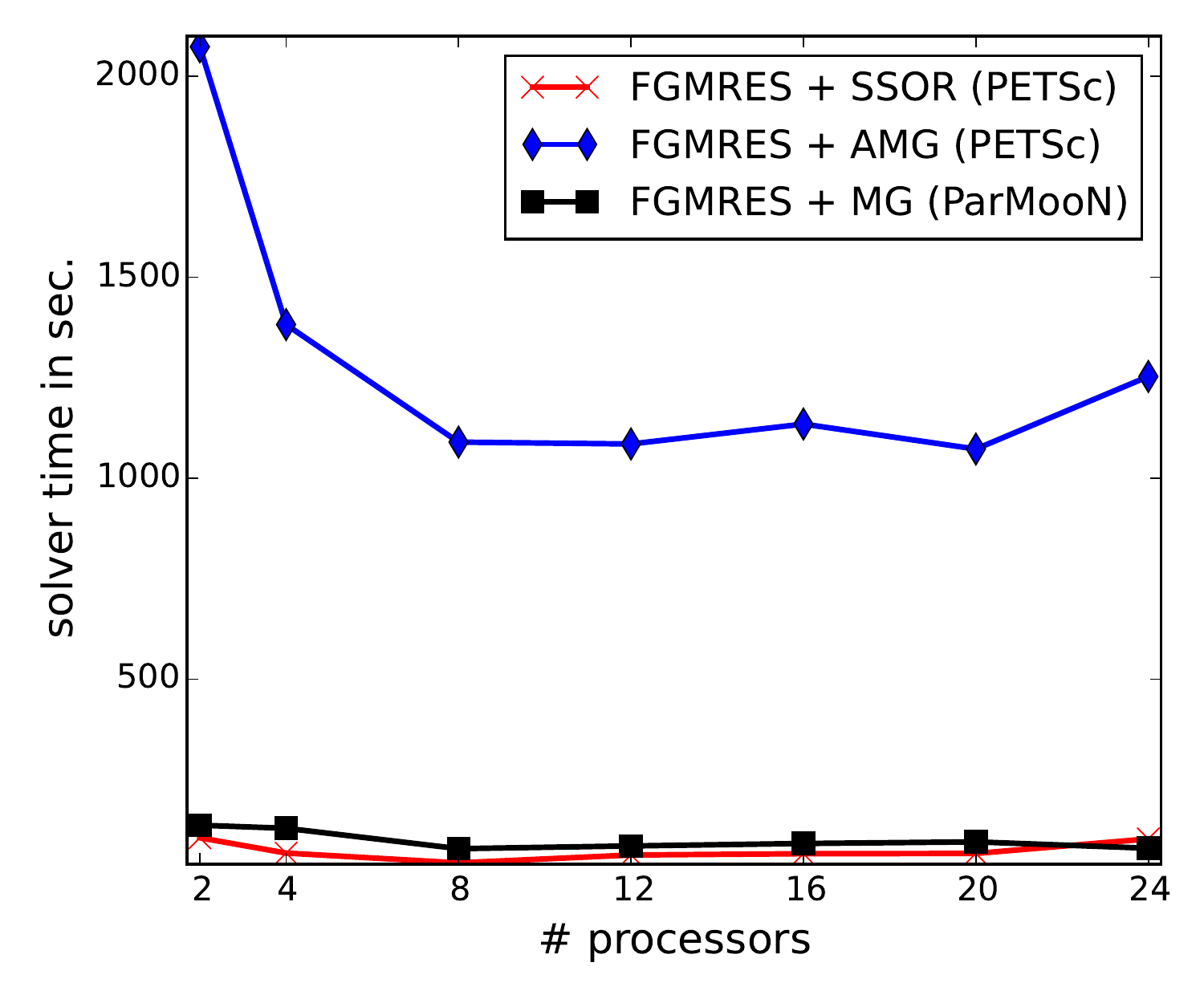}
\hspace*{1em} \includegraphics[width=0.48\textwidth]{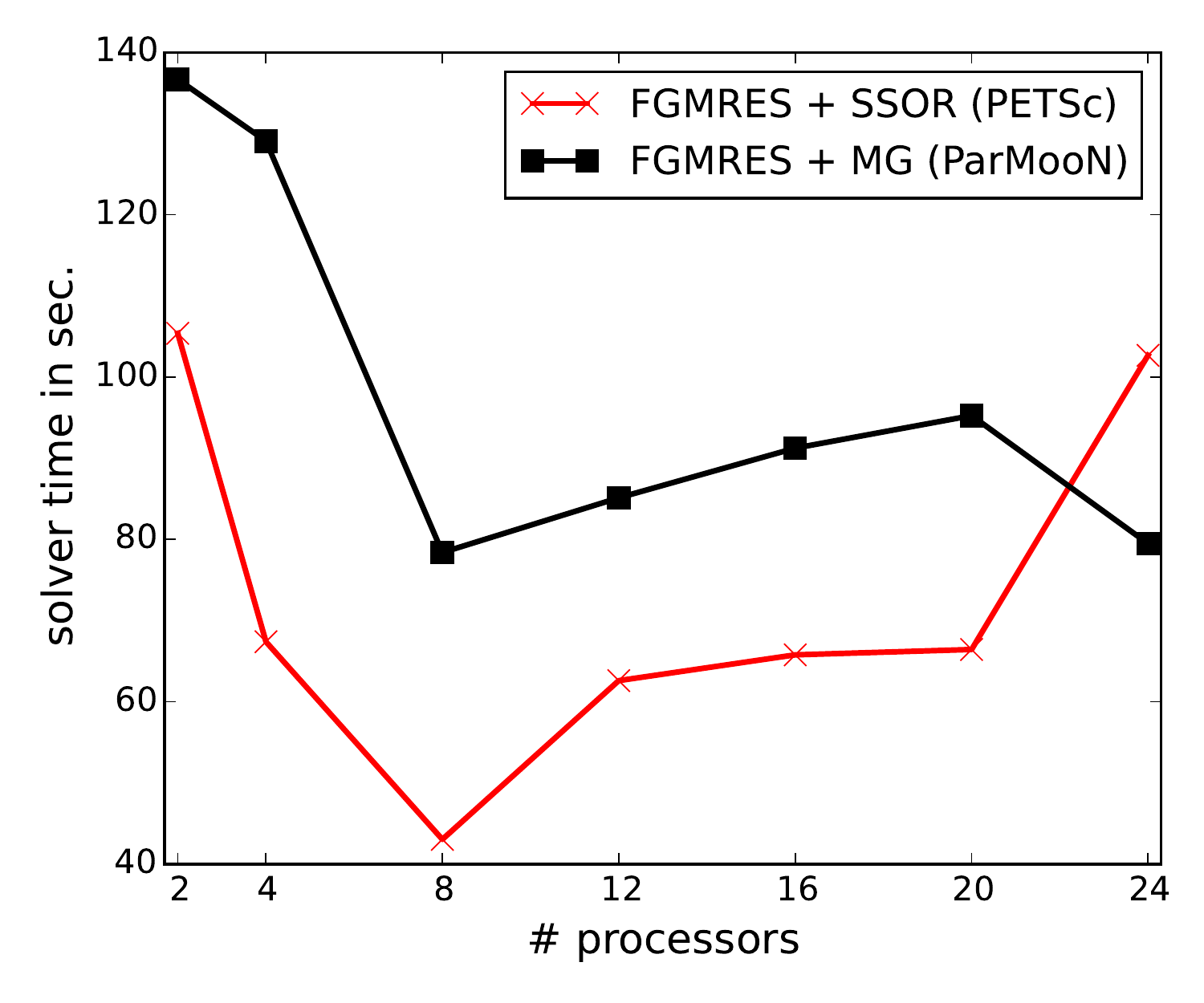}}
\caption{Example~\ref{ex:bail}, solver times on the $64^3$ cubed mesh.}
\label{fig:bail_results_64}
\end{figure}
\begin{figure}[t!]
\centerline{\includegraphics[width=0.48\textwidth]{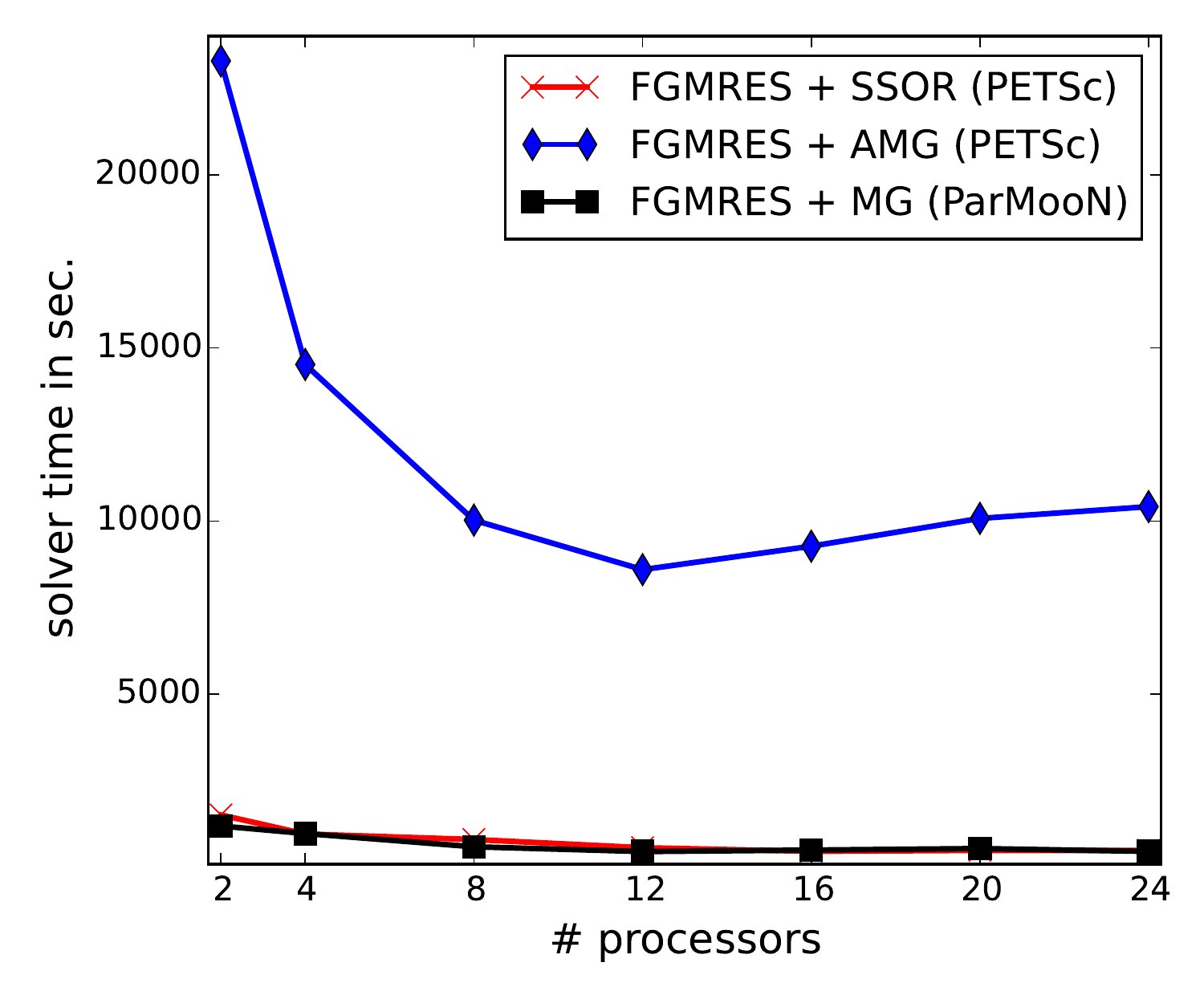}
\hspace*{1em} \includegraphics[width=0.48\textwidth]{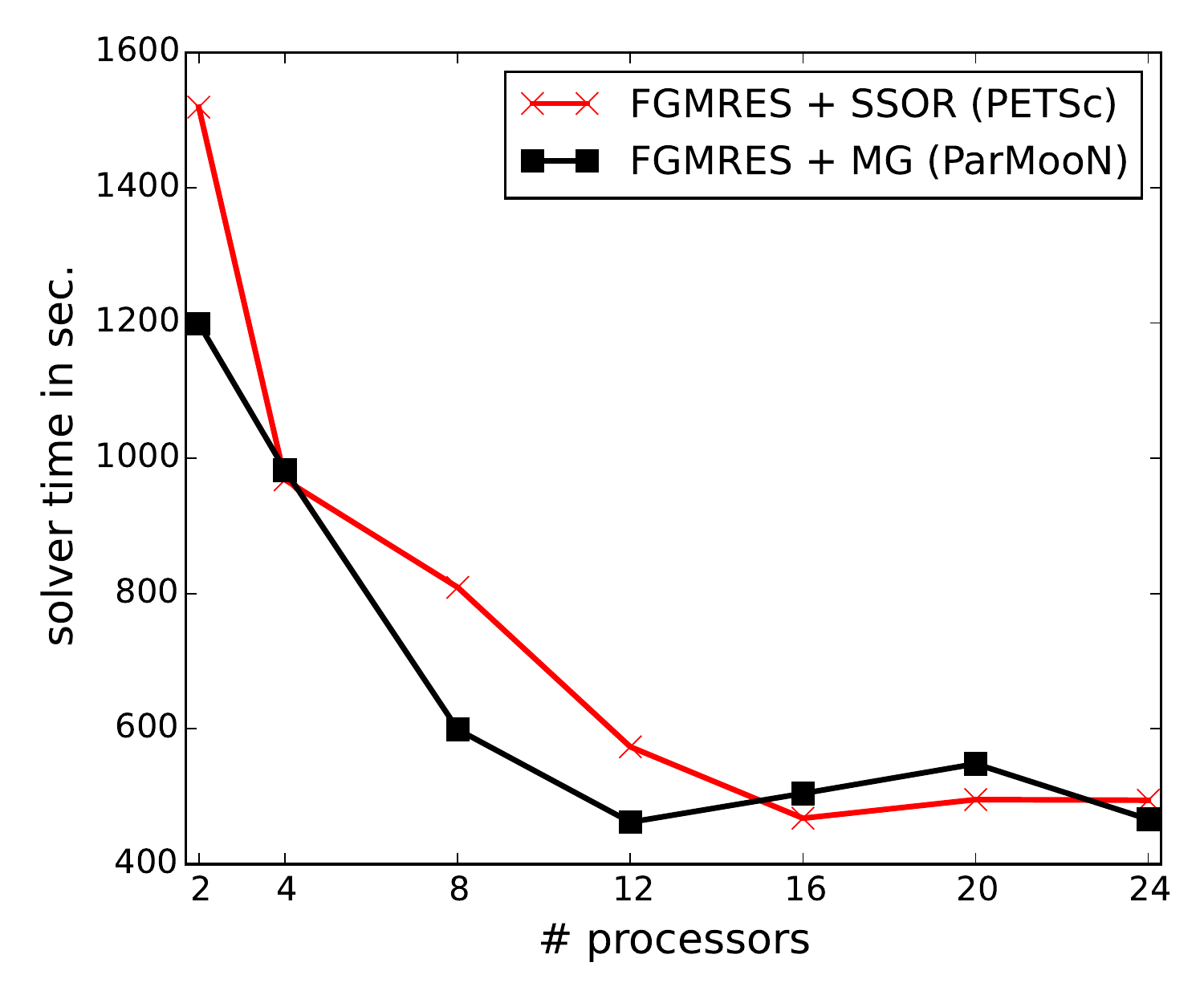}}
\caption{Example~\ref{ex:bail}, solver times on the  $128^3$ cubed mesh.}
\label{fig:bail_results_128}
\end{figure}
\end{example}

\begin{figure}[t!]
\centerline{\includegraphics[width=0.48\textwidth]{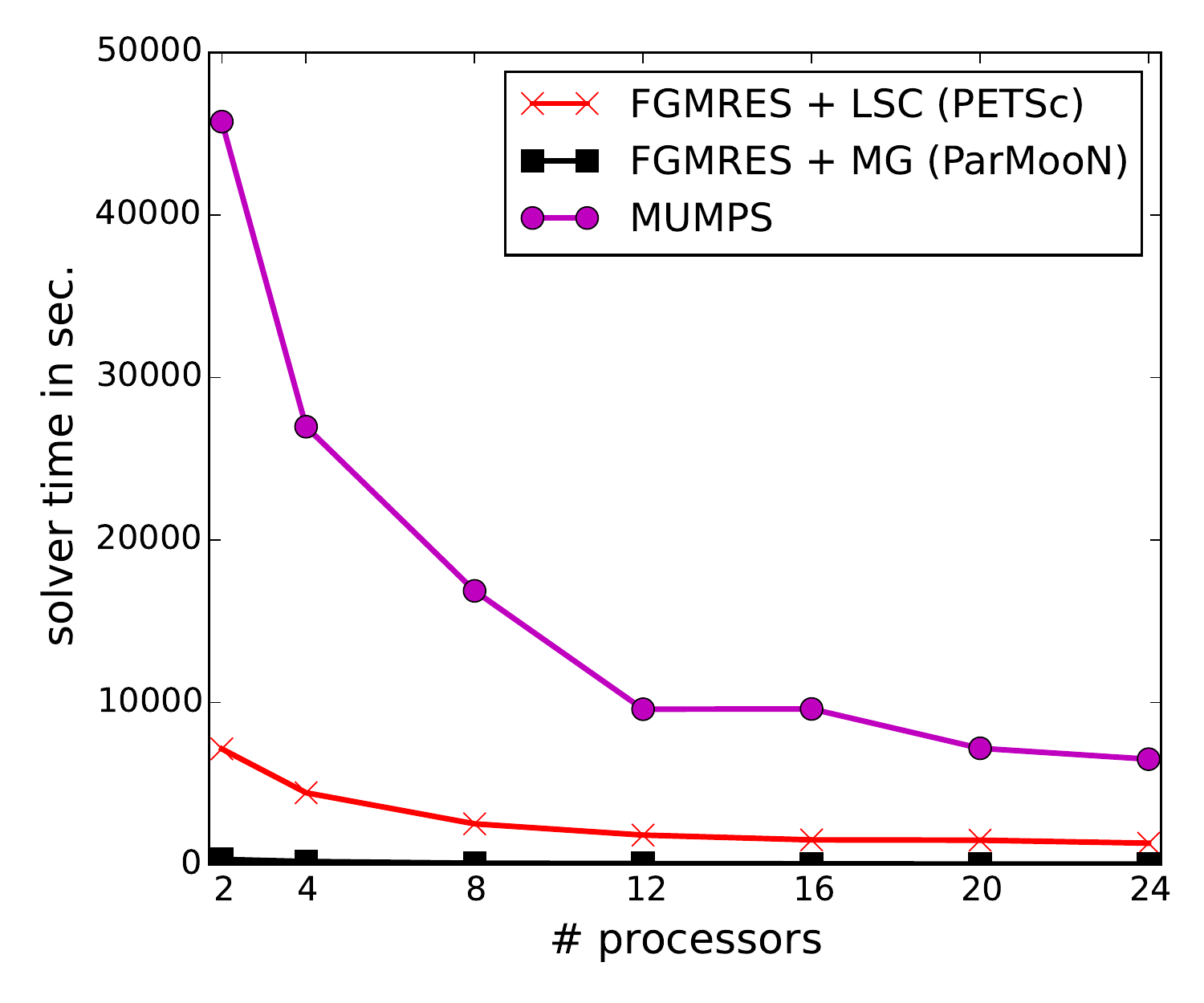}
\hspace*{1em} \includegraphics[width=0.48\textwidth]{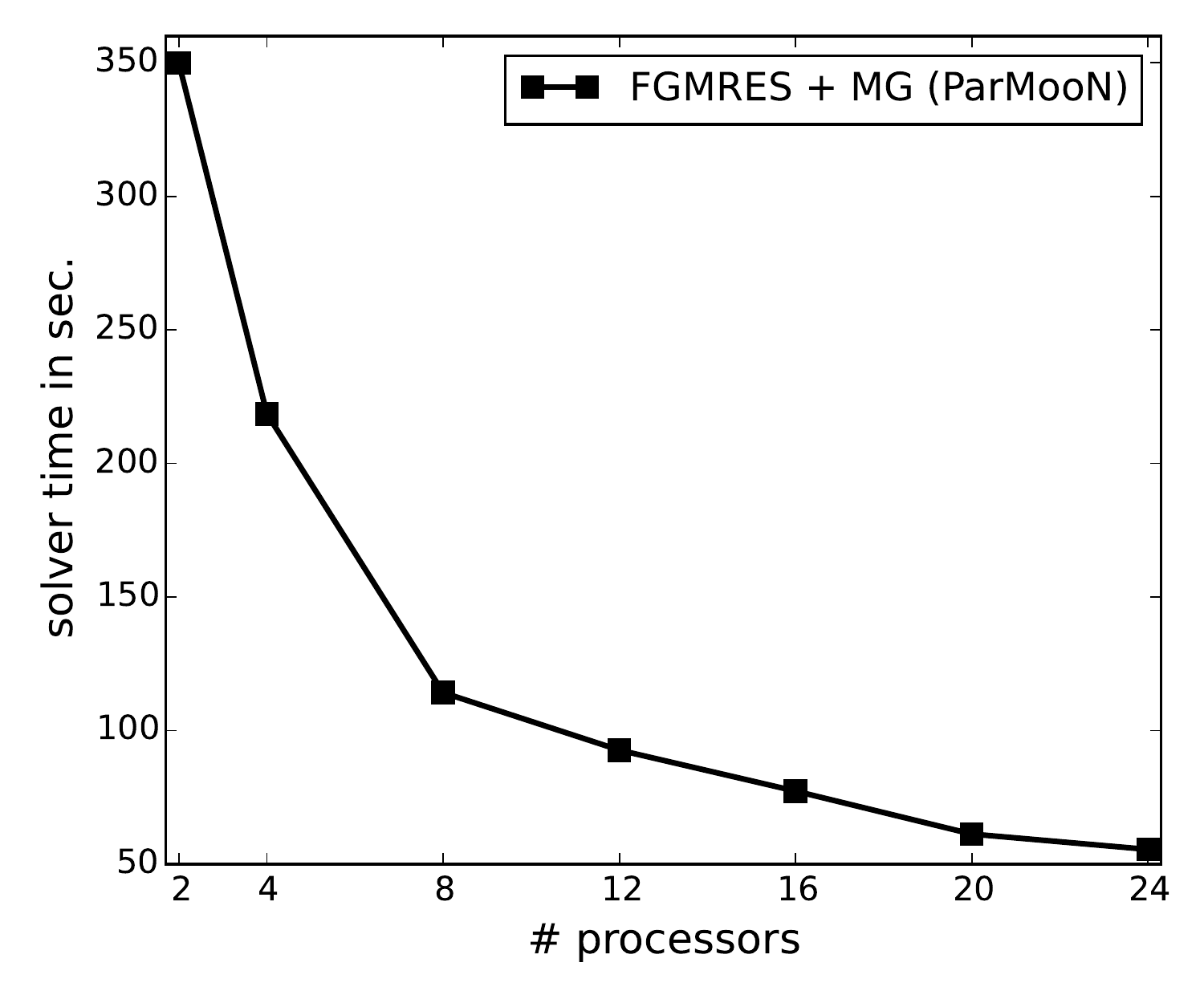}}
\caption{Example~\ref{ex:nse}, solver times on refinement level~2.}
\label{fig:nse_results_2}
\end{figure}

\begin{figure}[t!]
\centerline{\includegraphics[width=0.48\textwidth]{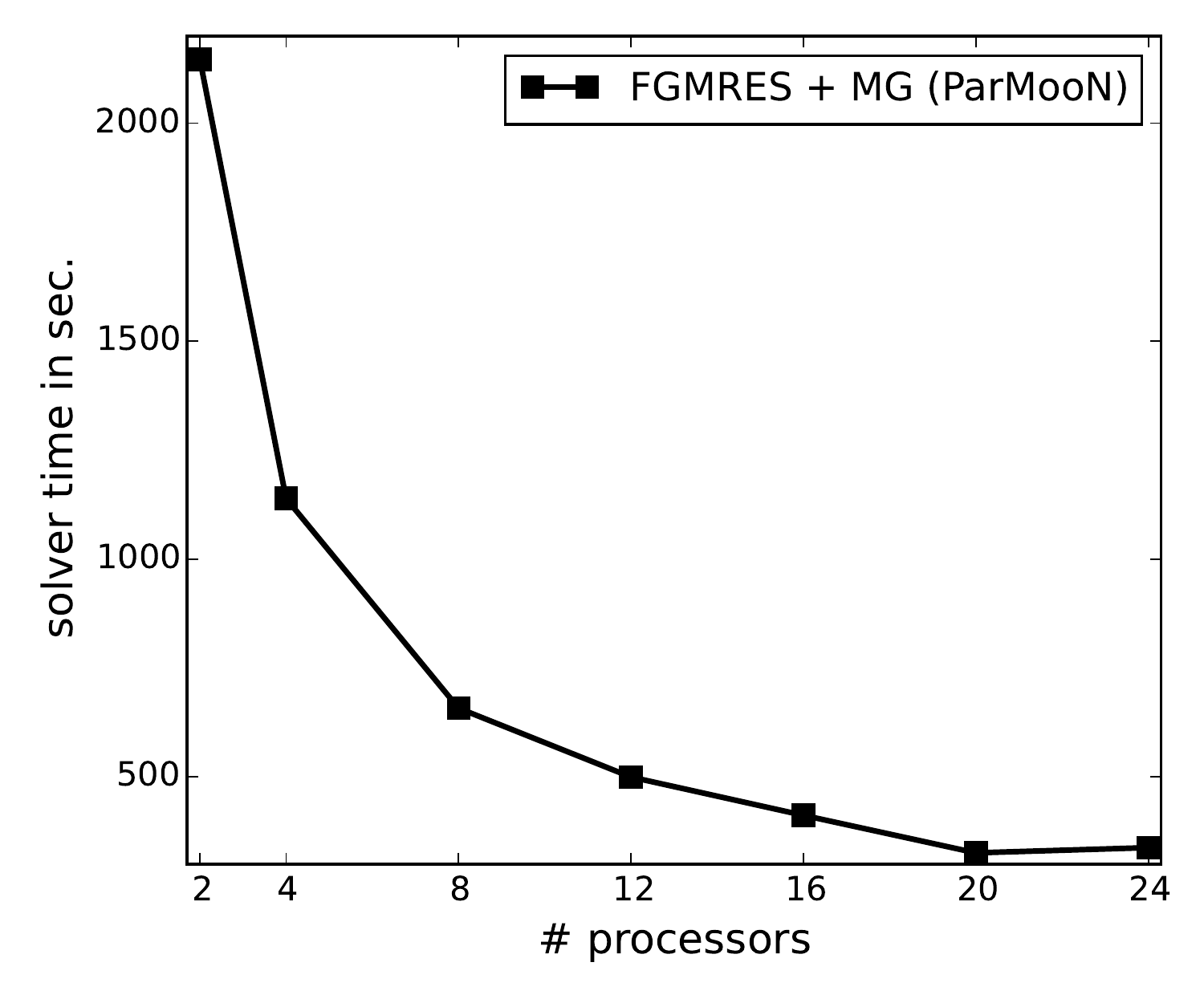}
\hspace*{1em} \includegraphics[width=0.48\textwidth]{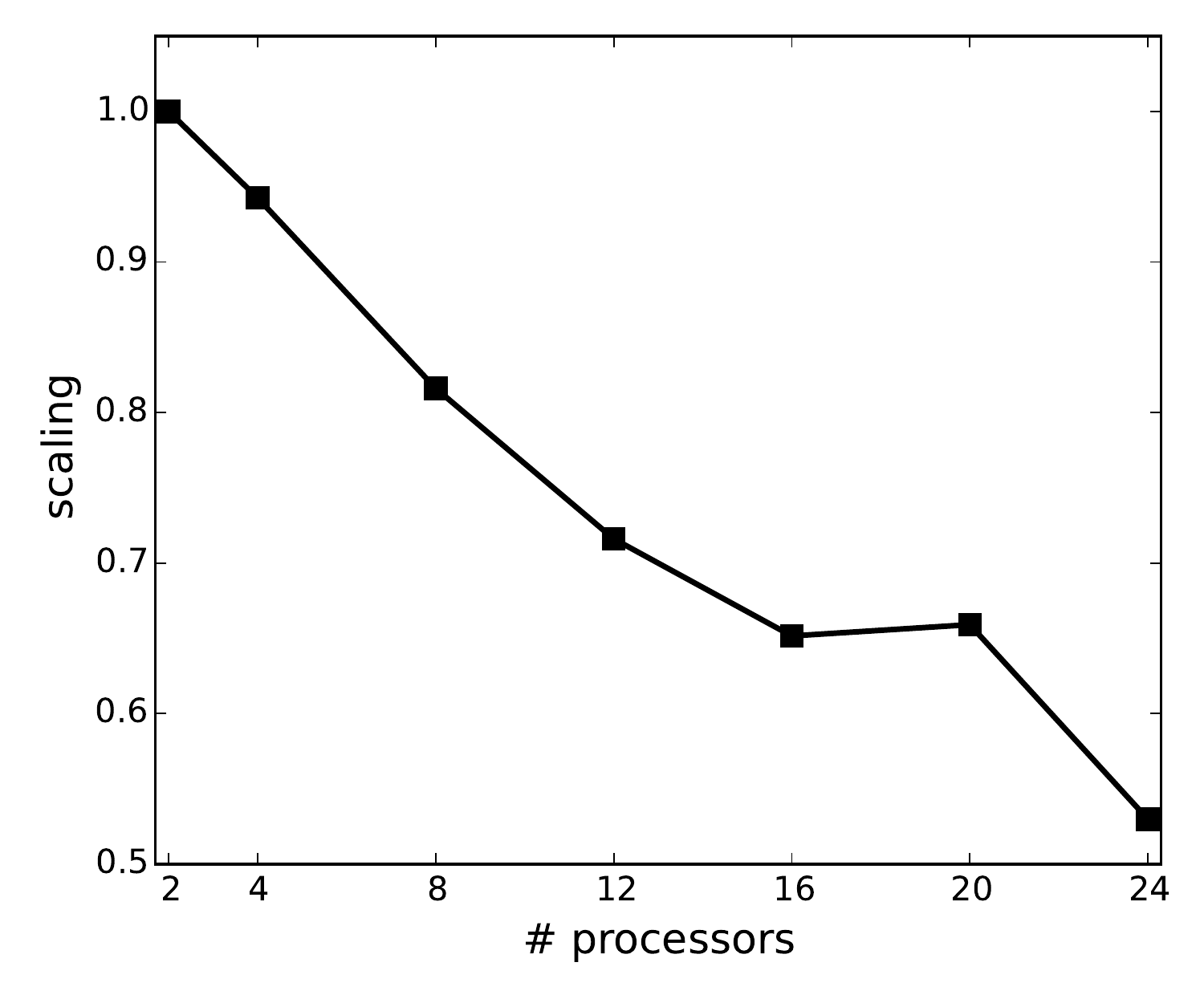}}
\caption{Example~\ref{ex:nse}, solver times and scaling on refinement level~3. The scaling is computed by 
$2\cdot t_2/(p\cdot t_p)$, where $p$ is the number of processors and $t_p$ the corresponding 
time from the left picture.}
\label{fig:nse_results_3}
\end{figure}

\begin{example}\rm {\em Steady-state incompressible Navier--Stokes equations.} 
\label{ex:nse} This example considers the benchmark problem of the flow around a cylinder defined in \cite{ST96}. 
The steady-state Navier--Stokes equations are given by 
$$
\begin{array}{rcll}
-\nu \Delta\bu + (\bu\cdot\nabla)\bu + \nabla p & = & \boldsymbol 0& \mbox{ in } \Omega,\\
\nabla\cdot\bu & = & 0 & \mbox{ in } \Omega,
\end{array}
$$
where $\bu$ is the velocity, $p$ the pressure, $\nu=10^{-3}$ is the dimensionless 
viscosity, and $\Omega$ is the domain given by
$$
\Omega = \left\{(0,2.55)\times (0,0.41) \setminus \left\{ (x-0.5)^2 + (y-0.2)^2 \le 0.05^2\right\} \right\}\times (0,0.41).
$$
At the inlet $x=0$, the velocity was prescribed by 
$$
\bu = \left( \frac{7.2}{0.41^4} yz(0.41-y)(0.41-z), 0 , 0 \right)^T,
$$
at the outlet $x=2.55$, the do-nothing boundary condition was imposed and on all other boundaries, 
the no-slip boundary condition was used. The flow field exhibits vortices behind the 
cylinder. 

The Navier--Stokes equations were discretized with the popular $Q_2/P_1^{\mathrm{disc}}$ 
(continuous piecewise triquadratic velocity, discontinuous piecewise linear pressure)
pair of finite element spaces on a hexahedral grid, see Figure~\ref{fig:nse_inigrid} 
for the initial grid. Simulations were performed on level~2 (776~160 velocity d.o.f.s, 122~800 
pressure d.o.f.s) and level~3 (6~052~800 velocity d.o.f.s, 983~040
pressure d.o.f.s). The nonlinear problem was linearized with a Picard iteration (fixed point 
iteration). In each step of this iteration, the iterative solvers reduced the Euclidean
norm of the residual vector by at least the factor 10 before performing the next Picard 
iteration. The Picard iteration was stopped if the Euclidean norm of the residual vector 
was less than $10^{-8}$.

\begin{figure}[t!]
\centerline{\includegraphics[width=0.67\textwidth]{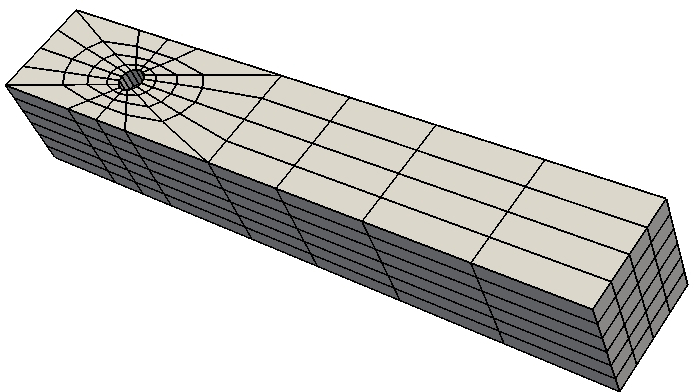}}
\caption{Example~\ref{ex:nse}, initial grid (level~0).}
\label{fig:nse_inigrid}
\end{figure}

In the geometric multigrid preconditioner, the so-called mesh-cell oriented Vanka 
smoother was used, e.g., see \cite{Joh02,Joh16}. This smoother is a block Gauss--Seidel method that 
solves a local system in each mesh cell. The multigrid preconditioner was applied with the 
V(2,2)-cycle. 

The linearization and used discretization of the incompressible Navier--Stokes equations require
the solution of a linear saddle point problem in each Picard iteration. We tried several options 
for calling an iterative solver from {\sc PETSc} for such problems, as well based on the coupled system 
as on Schur complement approaches. Only with a Schur complement approach and 
the least square commutator (LSC) preconditioner
({\tt -ksp\_type fgmres -pc\_type fieldsplit -pc\_fieldsplit\_type schur -pc\_fieldsplit\_schur\_factorization\_type upper -fieldsplit\_1\_pc\_type lsc 
-fieldsplit\_1\_lsc\_pc\_type lu -fieldsplit\_0\_ksp\_type preonly -fieldsplit\_0\_pc\_type lu}), 
reasonable computing times could be obtained, at least for the coarser grid, 
see Figure~\ref{fig:nse_results_2}. Like for the scalar problems, the sparse direct solver 
performed by far the least efficient. The solver with the geometric multigrid method 
was more efficient by around a factor of $20$ compared with the iterative solver 
called from {\sc PETSc}. For the finer grid, the numerical studies were restricted to the 
geometric multigrid preconditioner, see Figure~\ref{fig:nse_results_3}.
It can be observed
that increasing the number of processors from 2 to 20 reduces the computing time by 
a factor of around $6.6$.
\end{example}

\section{Summary}

This paper presented some aspects of the remanufacturing of an existing research code, in 
particular with respect to its parallelization. All distinct features of the predecessor
code could be incorporated in a straightforward way into the modernized code {\sc ParMooN}.
The efficiency of the most complex method in the parallel implementation, the geometric 
multigrid method, was assessed against some parallel solvers that are available in external 
libraries. The major conclusions of this assessment are twofold. For scalar convection-diffusion-reaction
equations, the geometric multigrid preconditioner was similarly efficient as an iterative solver 
from the {\sc PETSc} library. The larger the problems became, the better was its efficiency
in comparison with the external solver. For linear saddle point problems, arising in the 
simulation of the incompressible Navier--Stokes equations, we could not
find so far any external solver that proved to be nearly as efficient as the geometric
multigrid preconditioner. 

On the one hand, we keep on working at improving the efficiency of the geometric multigrid
preconditioner. On the other hand, we continue to assess external libraries with respect 
to efficient solvers for linear saddle point problems, which can be used in situations
where a multigrid hierarchy is not available. In addition, an assessment as presented
in this paper on another widely available hardware platform, namely small clusters of 
processors, is planned.

\section*{References}

\bibliographystyle{plain}
\bibliography{./main}

\end{document}